\documentclass[11pt]{amsart}

\pdfoutput=1

\usepackage{amsmath,amssymb,amsthm}
\usepackage{mathtools,url}
\usepackage{stmaryrd,xcolor,tikz-cd,textcomp}
\usepackage{enumitem}
\usepackage[top=1.25in, bottom=1in, left=1.25in, right=1.25in]{geometry}
\newcommand{\into}{\hookrightarrow}

\DeclareMathOperator{\Hom}{Hom}
\DeclareMathOperator{\Ext}{Ext}

\DeclareMathOperator{\Tor}{Tor}

\DeclareMathOperator{\Spec}{Spec}

\DeclareMathOperator{\Sym}{Sym}

\DeclareMathOperator{\Aut}{Aut}
\DeclareMathOperator{\ch}{ch}
\DeclareMathOperator{\orbch}{\wtilde{\ch}}
\DeclareMathOperator{\tr}{tr}
\DeclareMathOperator{\td}{td}
\DeclareMathOperator{\orbtd}{\wtilde{\td}}
\DeclareMathOperator{\ord}{ord}
\DeclareMathOperator{\Supp}{Supp}

\DeclareMathAlphabet{\mathpzc}{OT1}{pzc}{m}{it}

\newcommand{\vphi}{\varphi}
\newcommand{\wtilde}{\widetilde}

\newcommand{\inprod}[1]{\langle#1\rangle}
\newcommand{\tb}{\textbf}
\newcommand{\orbv}{\wtilde{v}}
\newcommand{\orbw}{\wtilde{w}}

\newcommand{\orbs}{\wtilde{s}}
\newcommand{\orbt}{\wtilde{t}}

\newcommand{\cj}{\overline}

\newcommand{\id}{\mathrm{id}}

\newcommand{\rk}{\mathrm{rk}}

\newcommand{\Coh}{\mathrm{Coh}}

\newcommand{\eig}{\mathrm{eig}}
\newcommand{\GL}{\mathrm{GL}}

\newcommand{\pt}{\mathrm{pt}}
\newcommand{\fix}{\mathrm{fix}}
\newcommand{\mov}{\mathrm{mov}}
\newcommand{\Tot}{\mathrm{Tot}}
\newcommand{\Ad}{\mathrm{Ad}}

\newcommand{\Sch}{\mathrm{Sch}}

\newcommand{\fppf}{\mathrm{fppf}}
\newcommand{\Ab}{\mathrm{Ab}}

\newtheorem{thm}{Theorem}[section]

\newtheorem{prop}[thm]{Proposition}
\newtheorem{lem}[thm]{Lemma}

\theoremstyle{definition}
\newtheorem{defn}[thm]{Definition}
\newtheorem{defn-prop}[thm]{Definition/Proposition}
\newtheorem{ex}[thm]{Example}

\newtheorem{rmk}[thm]{Remark}
\newtheorem{notation}[thm]{Notation}

\newcommand{\CC}{\mathbb{C}}

\newcommand{\PP}{\mathbb{P}}
\newcommand{\QQ}{\mathbb{Q}}

\newcommand{\ZZ}{\mathbb{Z}}

\newcommand{\cE}{\mathcal{E}}
\newcommand{\cF}{\mathcal{F}}

\newcommand{\cL}{\mathcal{L}}

\newcommand{\cO}{\mathcal{O}}

\newcommand{\cV}{\mathcal{V}}
\newcommand{\cW}{\mathcal{W}}
\newcommand{\cX}{\mathcal{X}}
\newcommand{\cY}{\mathcal{Y}}

\newcommand{\fg}{\mathfrak{g}}

\begin{document}
\title{Orbifold Hirzebruch-Riemann-Roch}
\author{Yuhang Chen}

\begin{abstract}
	We study the orbifold Hirzebruch-Riemann-Roch (HRR) theorem for quotient Deligne-Mumford stacks, explore its relation with the representation theory of finite groups, and derive a new orbifold HRR formula via an orbifold Mukai pairing.
\end{abstract}

\maketitle
% Show main sections only
\setcounter{tocdepth}{1}
\tableofcontents

\section{Introduction}
\vspace{1pt}
The Hirzebruch-Riemann-Roch (HRR) theorem is a fundamental result in algebraic geometry. Let $ X $ be a proper smooth scheme over $ \CC $. Let $ K(X) $ denote the Grothendieck ring of coherent sheaves on $ X, $ and let $ A(X) $ denote the Chow ring of $ X. $ Then we have a Chern character map $ \ch: K(X) \to A(X) \otimes {\QQ} $, which is a ring map, and a Todd class $ \td_X $, which is a unit in $ A(X) \otimes {\QQ} $. The HRR theorem says
\begin{equation}\label{HRR}
	\chi(X, x) = \int_X \ch(x) \td_X
\end{equation}
for any element $ x $ in $ K(X). $

The HRR theorem has been generalized to orbifolds by Kawasaki in \cite{kawasaki1979riemann} using analytic tools and later to Deligne-Mumford (DM) stacks by To\"en in \cite{toen1999theoremes} using algebraic methods. In both cases, we will call such a theorem an orbifold HRR theorem. The key difference between the orbifold HRR theorem and the ordinary one is that on the right side of formula (\ref{HRR}), the integration is carried out over the inertia stack $ I\cX $ of $ \cX $. Therefore, a good understanding of the inertia stack $ I\cX $ is crucial for one to apply the orbifold HRR theorem in practice. However, the treatment of $ I\cX $ in To\"en's paper is abstract, which makes applications of the oribifold HRR theorem difficult. In To\"en's work, it's assumed that a DM stack $ \cX $ has the resolution property, which implies that $ \cX $ is a quotient stack. For such a stack $ \cX $, Edidin \cite{edidin2013riemann} gave an explicit description of the inertia stack $ I\cX $, which is helpful for people working in the field to carry out explicit computations.
Let $ \cX $ be a connected proper smooth quotient DM stack over $ \CC $. Let $ K(\cX) $ be the Grothendieck ring of coherent sheaves on $ \cX, $ and let $ A(I\cX) $ be the Chow ring of $ I\cX $. Denote $ A(I\cX) \otimes \CC $ by $ A(I\cX)_\CC. $ Then we have an orbifold Chern character map $ \orbch: K(\cX) \to A(I\cX)_{\CC} $, which is a ring map, and an orbifold Todd class $ \orbtd_{\cX} $, which is a unit in $ A(I\cX)_{\CC}. $
The orbifold HRR theorem for $ \cX $ says
\begin{equation}\label{orb_HRR_thm}
	\chi(\cX, x) = \int_{I\cX} \orbch(x) \orbtd_{\cX}
\end{equation}
for any element $ x $ in $ K(\cX). $

In \cite{toen1999theoremes}, following the philosophy of Grothendieck, To\"en also proved a Grothendieck-Riemann-Roch theorem, which generalizes formula (\ref{orb_HRR_thm}) to a morphism of DM stacks. 
Instead of considering a morphism of DM stacks, we will generalize formula (\ref{orb_HRR_thm}) to a pair of coherent sheaves on $ \cX. $ We will introduce two new notions: an orbifold Mukai vector $ \orbv: K(\cX) \to A(I\cX)_{\CC} $, which is a linear map, and an orbifold Mukai pairing 
$$ \inprod{\cdot\ {,}\ \cdot}_{I\cX}: A(I\cX)_\CC \times A(I\cX)_\CC \to \CC, $$
which is sesquilinear, i.e., conjugate linear in the first component and linear in the second.
Let $ \chi(\cdot\ {,}\ \cdot) $ denote the orbifold Euler pairing on $ K(\cX) $.
We will prove that
\begin{equation*}\label{eq_orbifold_HRR_formula}
	\chi(x, y) = \inprod{\orbv(x),\orbv(y)}_{I\cX}
\end{equation*}
for any two elements $ x $ and $ y $ in $ K(\cX). $

\subsection{Outline}
The paper is structured as follows. In Section \ref{K-theory}, we review some basic facts about the algebraic $ K $-theory of quotient stacks. We will define the $ K $-theoretic Euler class of a vector bundle on a quotient stack $ \cX $ via the $ \lambda $-ring structure on the Gothendieck ring $ K^0(\cX) $ of vector bundles on $ \cX $. We will use the $ K $-theoretic Euler class to compute the Gothendieck ring of a weighted projective stack.

In Section \ref{Inertia_stacks}, we will describe a decomposition of the inertia stack $ I\cX $ of a connected separated quotient DM stack $ \cX = [X/G] $ indexed by a finite set of elements in $ G $. Each connected component in $ I\cX $ is a closed substack of $ \cX $ of the form $ \{g_i\} \times [X_i/G_i] $, where $ X_i $ is a closed subscheme of $ X $ invariant under the action of a subgroup $ G_i $ of the centralizer $ Z_{g_i} $ and fixed by an element $ g_i \in G_i $. We emphasize that the index $ g_i $ is important, as we will use it to define the orbifold Chern character and the orbifold Todd class map on $ K(\cX) $.

In Section \ref{Orbifold_Chern}, for a connected separated smooth quotient DM stack $ \cX, $ we will give detailed constructions for the orbifold Chern character map $ \orbch: K(\cX) \to A(I\cX)_\CC $, which is a $ \CC $-algebra isomorphism after tensored with $ \CC, $ and the orbifold Todd class map $ \orbtd: K(\cX) \to A(I\cX)_\CC^\times, $ which is a multiplicative map into the units in $ A(I\cX)_\CC. $ We will show that the orbifold Chern character map reduces to the inverse discrete Fourier transform when $ \cX $ is the classifying stack $ B\mu_n $ of the cyclic group $ \mu_n $.

In Section \ref{Orbifold_Mukai}, for a connected separated (resp. proper) smooth quotient DM stack $ \cX $, we will define the orbifold Mukai vector map (resp. the orbifold Mukai pairing). We will derive a new orbifold HRR formula for a pair of coherent sheaves on $ \cX. $ We will show that the orbifold HRR formula reduces to Parseval's theorem for the discrete Fourier transform when $ \cX = B\mu_n $.

\subsection{Conventions}
We fix the base field to be $ \CC $ throughout this paper. Let $ \pt = \Spec \CC. $ 

For a locally noetherian scheme $ X $, a coherent sheaf on $ X $ is a sheaf of $ \cO_X $-modules of finite presentation, and a vector bundle on $ X $ is a locally free sheaf of finite constant rank on $ X $. 

Let $ (\Sch /\CC)_\fppf $ denote the big fppf site of $ \Spec \CC $. An algebraic space is a sheaf $ X $ of sets on $ (\Sch /\CC)_\fppf $ such that the diagonal morphism $ X \to X \times X $ is representable by schemes, and there is a scheme $ U $ and a surjective \'etale morphism $ U \to X. $  A stack is a stack of groupoids over $ (\Sch /\CC)_\fppf $. An algebraic (resp. DM) stack is a stack $ \cX $ such that the diagonal morphism $ \cX \to \cX \times \cX $ is representable by algebraic spaces, and there is a scheme $ U $ and a surjective smooth (resp. \'etale) morphism $ U \to \cX. $

A quotient stack is a stack of the form $ [X/G] $, where $ X $ is a separated scheme of finite type, and $ G $ is a affine smooth group scheme of finite type (i.e., a linear algebraic group) acting on $ X $ from the left. 

For a group $ G, $ let $ 1 $ denote its identity. A representation of $ G $ is a linear representation of $ G. $ Let $ \mu_n $ denote the group scheme $ \Spec \CC[x]/(x^n-1). $ For a ring $ R, $ let $ R_\CC = R \otimes \CC. $

\subsection{Acknowledgements}
I would like to thank Promit Kundu and Xiping Zhang for many fruitful discussions. I want to thank Yonghong Huang for discussions on stacks and Zongzhu Lin for discussions on group schemes.
I also want to thank Dan Edidin and Johan de Jong for useful email correspondence.

\vskip 2em

\section{$ K $-theory of Quotient Stacks}\label{K-theory}
In this section we review some basic facts about the algebraic $ K $-theory of quotient stacks.
\vskip 4pt
Let $ \cX = [X/G] $ be a quotient stack. Then $ \cX $ is an algebraic stack with a smooth cover $ p: X \to \cX $, and $ \cX $ is a DM stack if the stabilizer of every geometric point of $ X $ is finite (and reduced, which is automatic for any group scheme over $ \CC $). Let $ \cX_\fppf $ denote the big fppf site of $ \cX $. Then we have a ringed site $ (\Ab(\cX_\fppf), \cO_\cX) $, where $ \Ab(\cX_\fppf) $ is the category of sheaves of abelian groups on $ \cX_\fppf $, and $ \cO_\cX $ is the structure sheaf of $ \cX $ which is a ring object in $ \Ab(\cX_\fppf) $. Now we can define coherent sheaves, locally free sheaves, and vector bundles on $ \cX $.

\begin{defn}
	A sheaf of $ \cO_\cX $-modules is a module object over $ \cO_\cX $ in $ \Ab(\cX_\fppf) $. A coherent sheaf (resp. locally free, vector bundle) on $ \cX $ is a sheaf $ \cE $ of $ \cO_\cX $-modules on $ \cX $ such that $ p^*\cE $ is a coherent sheaf (resp. locally free, vector bundle) on $ X. $
\end{defn}

Let $ \Coh(\cX) $ denote the category of coherent sheaves on $ \cX, $ and let $ \Coh^G(X) $ denote the category of $ G $-equivariant coherent sheaves on $ X. $
Note that for a sheaf $ \cE $ on $ \cX, $ the pullback along $ p: X \to \cX $ not only gives a sheaf $ E = p^*\cE $ on $ \cX, $ but also equips $ E $ with a $ G $-equivariant structure. Indeed, the functor $ p^*: \Coh(\cX) \to \Coh^G(X) $ induces an equivalence of categories between $ \Coh(\cX) $ and $ \Coh^G(X). $ For the convenience of the reader, we now give the definition of a $ G $-equivariant sheaf on $ X. $

\begin{defn}\label{defn_equiv_sheaves}
	Let $ p_2: G \times X \to X $ be the projection to $ X, $ and let $ \sigma: G \times X \to X $ be the action of $ G $ on $ X. $ A \tb{$ G $-equivariant structure} of a sheaf $ E $ on $ X $ is an isomorphism
	\begin{equation*}
		\phi: p_2^*E \xrightarrow{\sim} \sigma^*E
	\end{equation*}
	of sheaves on $ G\times X $ which satisfies the \tb{cocycle condition:}
	\begin{equation*}
		(\mu \times \id_X)^* \phi = (\id_G \times \sigma)^* \phi \circ p_{23}^* \phi,
	\end{equation*}
	which is an identity of isomorphisms of sheaves on $ G \times G \times X $ in the diagram
	\vspace{5pt}
	\begin{equation}\label{eq_equiv}
		\begin{tikzcd}[column sep=1.5em,row sep=3em]
			p_3^*E \arrow[rr, "(\mu \times \id_X)^* \phi", "\simeq"{below}] 
			\arrow[dr, "p_{23}^* \phi"{below left}, "\simeq"{above, sloped}] & & s^*E \\
			& r^*E \arrow[ur, "(\id_G \times \sigma)^* \phi"{below right}, "\simeq"{above, sloped}] &,
		\end{tikzcd}\vspace{8pt}
	\end{equation}
	where $ p_3, r, s $ are three morphisms from $ G \times G \times X $ to $ X $ defined by
	\begin{align*}
		p_3:\quad (g,h,x) & \quad \longmapsto \quad x, \\
		r: 	\quad (g,h,x) & \quad \longmapsto \quad hx, \\
		s: 	\quad (g,h,x) & \quad \longmapsto \quad ghx,
	\end{align*}
	and $ \mu \times \id_X,p_{23},\id_G \times \sigma $ are three morphisms from $ G \times G \times X $ to $ G \times X $ defined by
	\begin{align*}
		\mu \times \id_X:		\quad (g,h,x) & \quad \longmapsto \quad (gh,x), \\
		p_{23}: 				\quad (g,h,x) & \quad \longmapsto \quad (h,x), \\
		\id_G \times \sigma: 	\quad (g,h,x) & \quad \longmapsto \quad (g,hx).
	\end{align*}
	A sheaf $ E $ on  $ X $ with a $ G $-equivariant structure $ \phi $ is called a \tb{$ G $-equivariant sheaf} on $ X $ and is denoted by $ (E,\phi) $.
\end{defn}
A $ G $-equivariant morphism of $ G $-equivariant sheaves is a morphism of sheaves which is compatible with $ G $-equivariant structures.
\begin{defn}
	A \tb{$ G $-equivariant morphism} $ f: (E,\phi) \to (F,\psi) $
	is a morphism $ f: E \to F $ which is $ G $-equivariant, i.e., $ \psi \circ p_2^*f = \sigma^*f \circ \phi $ in the following diagram:
	\vspace{5pt}
	\begin{equation*}\label{eq_equiv_morphism}
		\begin{tikzcd}[column sep=3em,row sep=3em]
			p_2^*E \arrow{r}{p_2^*f} \arrow[d, "\phi"{left, outer sep = 2pt}, "\simeq"{above, sloped}] 
			& p_2^*F \arrow[d, "\psi"{right, outer sep = 2pt}, "\simeq"{below, sloped}] \\
			\sigma^*E \arrow{r}{\sigma^*f} & \sigma^*F
		\end{tikzcd}\vspace{8pt}
	\end{equation*}
\end{defn}

\begin{rmk}[Interpretation of $ G $-equivariant structures on vector bundles]
	Consider a $ G $-equivariant vector bundle $ (V, \phi) $ on $ X $. Let the same letter $ V $ denote its total space $ \Spec \Sym V^\vee $. Then the $ G $-equivariant structure $ \phi $ is the same as an action 
	\begin{equation*}
		\phi: G \to \Aut(V)
	\end{equation*}
	such that the projection $ \pi: V \to X $ is $ G $-equivariant, i.e., $ \pi \circ \phi(g) = g \circ \pi.$
	A $ G $-equivariant morphism $ f: (V, \phi) \to (W, \psi) $ between two $ G $-equivariant vector bundles is the same as a $ G $-equivariant morphism $ f: V \to W $
	of total spaces of vector bundles, i.e., a commutative diagram
	\vspace{5pt}
	\begin{equation}\label{eq_equiv_morphism_Vbundles}
		\begin{tikzcd}[column sep=3.3em,row sep=3.3em]
			V \arrow{r}{f} \arrow[d, "\phi(g)"{left}, "\simeq"{above, sloped}] & W \arrow[d, "\psi(g)"{right}, "\simeq"{below, sloped}] \\
			V \arrow{r}{f} & W
		\end{tikzcd}\vspace{8pt}
	\end{equation}
	for all $ g $ in $ G. $ Note that we don't require the morphism $ f $ to preserve fibers in diagram (\ref{eq_equiv_morphism_Vbundles}).
\end{rmk}

From now on, we will not distinguish a sheaf $ \cE $ on $ \cX $ and its corresponding $ G $-equivariant sheaf $ (E, \phi) $ on $ X. $

\begin{notation}
	Let $ K(\cX) $ denote the \tb{Grothendieck group} of coherent sheaves on $ \cX $. Denote the class of a coherent sheaf $ \cE $ in $ K(\cX) $ by $ [\cE] $. Let $ K^0(\cX) $ denote the Grothendieck group of vector bundles on $ \cX $. Denote the class of a vector bundle $ \cV $ in $ K^0(\cX) $ by $ [\cV]^0 $.
\end{notation}

\begin{rmk}
	There is a natural map 
	$$ K^0(\cX) \to K(\cX), \quad [\cV]^0 \mapsto [\cV] $$ 
	which may not be an isomorphism unless $ \cX $ has the resolution property.
\end{rmk}

\begin{defn}
	We say $ \cX $ has the \tb{resolution property} if every coherent sheaf on $ \cX $ is a quotient of a vector bundle on $ \cX $.
\end{defn}
It's known that all factorial separated noetherian schemes have the resolution property. In particular, all separated smooth schemes of finite type have the resolution property. Since we assume X is separated of finite type in a quotient stack $ [X/G], $ the following holds.
\begin{prop}
	All smooth quotient stacks have the resolution property.
\end{prop}
\begin{proof}
	Let $ [X/G] $ be a smooth quotient stack. Then $ X $ is a separated smooth scheme of finite type, which ensures that $ [X/G] $ has the resolution property by \cite[Theorem 2.1]{totaro2004resolution}.
\end{proof}
\begin{rmk}
	The group $ K^0(\cX) $ always has a ring structure with the multiplication given by the tensor product of vector bundles and the unity $ 1 = [\cO_\cX]^0 $. We call $ K^0(\cX) $ the \tb{Grothendieck ring} of $ \cX. $
\end{rmk}
Now assume $ \cX $ is smooth. Then every coherent sheaf $ \cE $ on $ \cX $ has a finite resolution
\begin{equation}\label{eq_finite_resolution}
	0 \to \cE_n \to \cE_{n-1} \to \cdots \to \cE_1 \to \cE_0 \to \cE \to 0
\end{equation}
which is an exact sequence in $ \Coh(\cX) $, where each $ \cE_i $ is a vector bundle.
Since $ \cX $ is smooth, the natural map $ K^0(\cX) \to K(\cX) $ becomes an isomorphism of abelian groups with the inverse given by
\begin{equation*}
	\beta: K(\cX) \to K^0(\cX), \quad [\cE] \mapsto \sum_{i=0}^n (-1)^i[\cE_i]^0
\end{equation*}
for a finite resolution (\ref{eq_finite_resolution}) of a coherent sheaf $ \cE $. Therefore, we can equip the group $ K(\cX) $ with a ring structure. We also call $ K(\cX) $ the Grothendieck ring of $ \cX. $ Alternatively, we can define a multiplication directly on $ K(\cX) $ as follows.
\begin{defn}\label{defn_multiplication_sheaves}
	Let $ \cX $ be a smooth quotient stack. Define a multiplication on $ K(\cX) $ by
	\begin{equation}\label{eq_multiplication_sheaves}
		[\cE] [\cF] = \sum_i (-1)^i [\Tor_i^{\cO_\cX}(\cE,\cF)]
	\end{equation}
	for coherent sheaves $ \cE $ and $ \cF $ on $ \cX $ and extend linearly. 
\end{defn}
\begin{rmk}
	Note that the product of $ [\cE] $ and $ [\cF] $ above corresponds to the class of the derived tensor product $ \cE \otimes^L \cF $ in the Grothendieck group of the derived category of coherent sheaves on $ \cX. $
\end{rmk}
The following proposition is immediate.
\begin{prop}\label{lem_multiplication_sheaves}
	Let $ \cX  $ be a smooth quotient stack. Let $ \cE $ and $ \cF $ be coherent sheaves on $ \cX $. Then their product in $ K(\cX) $ is given by
	\begin{equation*}
		[\cE] [\cF] = \sum_i (-1)^i [\cE_i \otimes \cF]
	\end{equation*}
	for any finite resolution $ \cE_{\boldsymbol{\cdot}} \to \cE \to 0 $.
	If either $ \cE $ or $ \cF $ is locally free, then $ [\cE] [\cF] = [\cE \otimes \cF] $.
\end{prop}
The multiplication on $ K(\cX) $ defined by (\ref{eq_multiplication_sheaves}) is compatible with the multiplication on $ K^0(\cX) $ defined by the tensor product of vector bundles, and hence is
commutative and associative. 

From now on, we will omit the superscript $ ``0" $ in the notation $ [\cE]^0 $ for the class of a vector bundle $ \cE $ on $ \cX $. Now we drop the smoothness assumption of $ \cX. $

The Grothendieck ring $ K^0(\cX) $ has an additional algebraic structure called a $ \lambda $-ring. We follow \cite[Section I.1]{fulton1985riemann} for the definition of a $ \lambda $-ring, which is called a pre-$ \lambda $-ring in some literatures, for example in \cite{yau2010lambda}.
\begin{defn}
	A \tb{$ \lambda $-ring} is a ring $ R $ with a group homomorphism
	\begin{align*}
	\lambda_t: (R, +) &\to \left(R[[t]]^\times, \times\right)\\
	x & \mapsto \sum_{i \geq 0} \lambda^i(x) t^i
	\end{align*}
	such that $ \lambda^0(x) = 1 $ and $ \lambda^1(x) = x $ for all $ x \in K. $
	A \tb{positive structure} on a $ \lambda $-ring $ (R, \lambda_t) $ is a pair $ (\varepsilon, R_+) $, where $ \varepsilon: R \to \ZZ $ is a surjective ring homomorphism, and $ R_+ $ is a subset of $ R $ such that the following conditions hold:
	\begin{enumerate}[font=\normalfont,leftmargin=*]
		\item $ R_+ $ is closed under addition and multiplication.
		\item $ R_+ $ contains the set $ \ZZ_+ $ of positive integers.
		\item $ R = R_+ - R_+, $ i.e., every element of $ R $ is a difference of two elements of $ R_+ $.
		\item If $ x \in R_+ $, then $ \varepsilon(x) = r > 0,\ \lambda^i(x) = 0 $ for $ i > r $, and $ \lambda^r(x) $ is a unit in $ R. $
	\end{enumerate}
	Elements of $ R_+ $ are said to be \tb{positive}. An element of $ R $ is \tb{non-negative} if it is either positive or zero. The set of non-negative elements is denoted by $ R_{\geq 0}. $
\end{defn}
\begin{rmk}
	Consider a $ \lambda $-ring $ (R, \lambda_t). $ For all $ x, y \in R, $ we have
	\begin{equation*}
	\lambda_t(x+y) = \lambda_t(x)\lambda_t(y),
	\end{equation*}
	which is equivalent to
	\begin{equation*}
	\lambda^n(x+y) = \sum_{i=0}^{n} \lambda^i(x)\lambda^{n-i}(y)
	\end{equation*}
	for all integer $ n \geq 0.$ Note that we always have $ \lambda_t(0) = 1. $ Suppose $ (R, \lambda_t) $ has a positive structure $ (\varepsilon, R_+). $ Then we have $ \varepsilon(1) = 1 $ and $ \lambda_t(1) = 1 + t. $ On the other hand, we have $ \varepsilon(-1) = -1 $ and 
	\begin{equation*}
	\lambda_t(-1) = \lambda_t(1)^{-1} = 1 - t + t^2 - t^3 + \cdots.
	\end{equation*}
\end{rmk}
From now on, we always assume there is a positive structure associated with a $ \lambda $-ring.
\begin{ex}
	A simple example of a $ \lambda $-ring is $ \ZZ $ with $ \lambda_t(m) = (1+t)^m $ and the obvious positive structure $ (\varepsilon = \id: \ZZ \to \ZZ, \ZZ_+). $ In this case, we have 
	$$ \lambda^i(m) = \binom{m}{i} $$ 
	for $ m \in \ZZ $ and $ i \geq 0. $
\end{ex}

\begin{ex}[{\cite[Section V.1]{fulton1985riemann}}]
	Let $ X $ be a scheme. The Grothendieck ring $ K^0(X) $ of vector bundles on $ X $ is a $ \lambda $-ring: positive elements $ K^0(X)_+ $ are classes of vector bundles $ V $, taking exterior powers induces a map $ \lambda_t: K^0(X)_{\geq 0} \to K^0(X)[[t]]^\times $ given by
	\begin{equation*}
	\lambda_t[V] = \sum_{i \geq 0} \lambda^i [V] t^i = \sum_{i \geq 0} [\wedge^i V] t^i,
	\end{equation*}
	which is a semigroup homomorphism: for two vector bundles $ V $ and $ W,$
	\begin{equation*}
	\lambda_t([V] + [W]) = \lambda_t[V]\lambda_t[W]
	\end{equation*}
	because 
	\begin{equation*}
	\wedge^n(V \oplus W) \cong \bigoplus_{i=0}^{n} \left(\wedge^i V \otimes \wedge^{n-i} W\right) \ \text{for all integers} \ n \geq 0.
	\end{equation*}
	Since each element in $ K^0(X) $ can be written as a difference $ [V] - [W] $ of two classes of vector bundles, we can extend the domain of $ \lambda_t $ to all of $ K^0(X) $, i.e.,
	\begin{equation*}
	\lambda_t([V]-[W]) = \lambda_t[V]/\lambda_t[W].
	\end{equation*}
	The surjective ring homomorphism $ \varepsilon = \rk:  K^0(X) \to \ZZ $ is the rank map: for an element $ x = [V] - [W], $
	$$ \rk(x) = \rk(V) - \rk(W). $$ 
	When $ X = \pt, $ we get back the example of $ \ZZ. $ 
\end{ex}
\begin{prop}
	Let $ \cX = [X/G] $ be a quotient stack. The Grothendieck ring $ K^0(\cX) $ of vector bundles on $ \cX $ is a $ \lambda $-ring.
\end{prop}
\begin{proof}
	Consider a vector bundle $ \cV = (V, \phi) $ on $ \cX $. Exterior powers $ \wedge^i V $ of $ V $ carry canonical $ G $-equivariant structures $ \wedge^i \phi, $ so $ \wedge^i \cV $ is well-defined. Every element in $ K^0(\cX) $ is a difference $ [\cV] - [\cW] $ of two classes of vector bundles on $ \cX $, so we can define the same $ \lambda_t $ and $ \varepsilon = \rk $ as those in the case of schemes. 
\end{proof}
\begin{defn}
	Let $ \cX $ be a quotient stack. Taking dual vector bundles respects exact sequences and hence defines an involution 
	$$ (\ \cdot \ )^\vee: K^0(\cX) \to K^0(\cX). $$
\end{defn}
\begin{rmk}
	The involution $ (\ \cdot \ )^\vee: K^0(\cX) \to K^0(\cX) $ is a ring automorphism, i.e., for all $ x $ and $ y $ in $ K^0(\cX), $ we have
	\begin{equation*}
	1^\vee = 1, \quad (x+y)^\vee = x^\vee + y^\vee, \quad (xy)^\vee = x^\vee y^\vee, \quad \text{and} \quad x^{\vee\vee} = x.
	\end{equation*}
\end{rmk}
\begin{defn}
	Let $ \cX $ be a quotient stack. The \tb{$ K $-theoretic Euler class} of a non-negative element $ x = [\cV] $ in $ K^0(\cX) $ of rank $ r $ is defined by
	\begin{equation*}
	e^K(x) = \lambda_{t}(x^\vee) |_{t=-1} = \sum_{i=0}^{r} (-1)^i \lambda^i(x^\vee) = \sum_{i=0}^{r} (-1)^i [\wedge^i \cV].
	\end{equation*}
\end{defn}
\begin{rmk}[Properties of the $ K $-theoretic Euler class]
	Let $ \cX $ be a quotient stack. The $ K $-theoretic Euler class gives a semigroup homomorphism 
	$$ e^K: (K^0(\cX)_{\geq 0}, +) \to (K^0(\cX), \times), $$ 
	i.e., for all non-negative elements $ x $ and $ y $ in $ K^0(\cX), $ we have
	\begin{equation*}
	e^K(x+y) = e^K(x)e^K(y).
	\end{equation*}
	In particular, $ e^K(0) = 1. $ By definition, we also have $ e^K(1) = 1 - 1 = 0. $ 
	Note that $ e^K $ cannot be defined on all of $ K^0(\cX) $: for example, $ e^K(-1) $ is undefined since $ \lambda_t(-1) $ is an infinite series. If $ x $ and $ y $ are positive elements of $ K^0(\cX) $ such that $ e^K(y) $ is a unit, then we define
	\begin{equation*}
	e^K(x-y) = e^K(x)/e^K(y).
	\end{equation*}
\end{rmk}
\begin{rmk}[The geometric meaning of the $ K $-theoretic Euler class]
	Consider a vector bundle $ \cV = (V, \phi) $ of rank $ r $ on a quotient stack $ \cX = [X/G] $. Suppose the vector bundle $ \cV $ has a section $ s: \cO_\cX \to \cV $ (i.e., a $ G $-equivariant section $ s: \cO_X \to V $) such that it cuts out a substack $ \cY = [Y/G] $, where $ Y $ is a $ G $-invariant subscheme of $ X $ of codimension $ r. $ Let $ i: Y \to X $ denote the inclusion, and let $ s^\vee: V^\vee \to \cO_X $ denote the dual of $ s $. The Koszul complex
	\begin{equation*}
	0 \to \wedge^r V^\vee \to \wedge^{r-1} V^\vee \to \cdots \to V^\vee \xrightarrow{s^\vee} \cO_X \to i_*\cO_Y \to 0
	\end{equation*}
	on $ X $ is a $ G $-equivariant locally free resolution of the pushforward $ i_*\cO_Y, $ so we obtain a Koszul complex
	\begin{equation*}
	0 \to \wedge^r \cV^\vee \to \wedge^{r-1} \cV^\vee \to \cdots \to \cV^\vee \to \cO_\cX \to i_*\cO_\cY \to 0
	\end{equation*}
	on $ \cX. $
	Therefore,
	\begin{equation*}\label{eq_eK_V}
	[i_*\cO_\cY] = e^K(\cV),
	\end{equation*}
	which explains the name ``$ K $-theoretic Euler class" since the usual Euler class $ e(\cV) $ is the class 
	$$ [\cY] = c_r(\cV)$$
	in the $ r $-th Chow group $ A^r(\cX) $ of $ \cX. $
	% (or cohomology class of the Poincar\'e dual)
	Suppose $ \cX $ and $ \cY $ are smooth. Taking the equivariant Chern class as in Example 15.3.1 in \cite{fulton1998intersection}, we obtain a relation between the $ K $-theoretic Euler class and the usual Euler class:
	\begin{equation*}
	c_r(e^K(\cV)) = (-1)^{r-1} (r-1)! e(\cV).
	\end{equation*}
\end{rmk}
We can use the $ K $-theoretic Euler class to compute the Grothendieck ring of a weighted projective stack $ \PP(a_0, \dots, a_n) $.
\begin{ex}[The Grothendieck ring of a weighted projective stack]
	Let $ X = \CC^{n+1} $ with a closed subscheme $ p = \{0\} $, and let $ \cX = [X/\CC^*] $, where the $ \CC^* $-action on $ X $ is given by
	\begin{equation*}
	t (x_0, \dots, x_n) = (t^{a_0} x_0, \dots, t^{a_n} x_n)
	\end{equation*}
	for some positive integers $ a_0, \dots, a_n. $ Let $ U = X - p. $ We denote the weighted projective stack $ \PP(a_0, \dots, a_n) = [U/\CC^*] $ by $ \PP. $ Let $ i: p \into X $ and $ j: U \into X $ be the inclusions of $ p $ and $ U $ in $ X, $ which are $ G $-equivariant morphisms of schemes. Hence they induce inclusions 
	\begin{equation*}
	i: B\CC^* \to \cX \quad \text{and} \quad j: \PP \to \cX
	\end{equation*}
	of substacks in $ \cX. $
	Therefore, we have a right exact sequence
	\begin{equation}\label{eq_K_theory_right_exact}
	K(B\CC^*) \to K(\cX) \to K(\PP) \to 0
	\end{equation}
	of abelian groups, where the first map is extension by zero and the second one is restriction. Since $ \cX $ and $ \PP $ are smooth quotient stacks with the resolution property, we can identify the two maps in (\ref{eq_K_theory_right_exact}) as the $ K $-theoretic pushforward and pullback
	$$ i_K: K(B\CC^*) \to K(\cX) \quad \text{and} \quad j^K: K(\cX) \to K(\PP), $$ 
	where $ K(\cX) $ is identified with $ K^0(\cX). $
	We know that 
	\begin{equation*}
	K(B\CC^*) \cong \ZZ[u, u^{-1}] \quad \text{and} \quad K(\cX) \cong \ZZ[x, x^{-1}],
	\end{equation*}
	where $ u $ is the class of the identity character 
	$$ \rho = \id: \CC^* \to \CC^*, $$ 
	and $ x $ is the class of the $ \CC^* $-equivariant line bundle $ \cO_\cX \otimes \rho. $
	The tangent bundle $ T\cX \cong \cO_\cX \otimes (\rho^{a_0} \oplus \cdots \oplus \rho^{a_n}) $ on $ \cX $ has a canonical $ \CC^* $-equivariant section 
	$$ s: \cO_X \to TX, \quad 1 \mapsto (x_0, \dots, x_n), $$ 
	which cuts out the origin $ p $ of $ X. $
	Since
	\begin{equation*}
	[TX] = x^{a_0} + \cdots + x^{a_n},
	\end{equation*}
	the image of $ 1 = [\cO_{B\CC^*}] $ under the map $ i_K $ is
	\begin{equation*}
	i_K (1) = e^K(TX) = e^K \left(\sum_{i=0}^n x^{a_i}\right) = \prod_{i=0}^n e^K(x^{a_i}) = \prod_{i=0}^n \lambda_{-1}(x^{-a_i}) = \prod_{i=0}^n (1-x^{-a_i}).
	\end{equation*}
	Therefore,
	\begin{equation*}
	K(\PP) \cong \frac{\ZZ[x,x^{-1}]}{\inprod{(1-x^{-a_0})\cdots(1-x^{-a_n})}} \cong \frac{\ZZ[x]}{\inprod{(x^{a_0}-1)\cdots(x^{a_n}-1)}}.
	\end{equation*}
	Here $ x $ is interpreted as the class of the twisting sheaf
	$$ \cO_{\PP}(1) = j^*(\cO_\cX \otimes \rho) $$ 
	on $ \PP, $ and is invertible in $ K(\PP). $ Here each $ 1 - x^{-a_i} $ is the class of $ \cO_{\PP_i} $ in the short exact sequence
	\begin{equation*}
	0 \to \cO_\PP(-a_i) \to \cO_\PP \to \cO_{\PP_i} \to 0
	\end{equation*}
	of coherent sheaves on $ \PP $, where $ \PP_i = \PP(a_0, \dots, \hat{a_i},\dots,a_n) $ is the substack in $ \PP $ cut out by the coordinate $ x_i. $ Therefore, the relation
	\begin{equation*}\label{eq_wps_relation}
	(1-x^{-a_0}) \cdots (1-x^{-a_n}) = 0
	\end{equation*}
	in $ K^0(\PP) $ says
	\begin{equation*}
	[\cO_{\PP_0}] \cdots [\cO_{\PP_n}] = 0
	\end{equation*}
	in $ K^0(\PP) $, which reflects the fact that the intersection of $ \PP_0, \dots, \PP_n $ is empty.
\end{ex}
\vskip 2em

\section{Inertia Stacks of Quotient Stacks}\label{Inertia_stacks}
In this section we study the inertia stack associated to a quotient stack, which is the key ingredient for the orbifold HRR theorem.
\vskip 4pt
Let $ \cX = [X/G] $ be a quotient stack. There is a diagonal action of $ G $ on $ X \times X $ given by $ g(x,y) = (gx, gy) $ for $ g \in G $ and $ x, y \in X. $ 
Since $ G $ acts on itself by conjugation, there is also a diagonal action of $ G $ on $ G \times X $ given by $ h(g, x) = (hgh^{-1}, hx) $ for $ g, h \in G $ and $ x \in X. $
There is also an \tb{action morphism}, i.e., 
\begin{equation*}
	\alpha: G \times X \to X \times X, \quad (g,x) \mapsto (x, gx),
\end{equation*}
which is $ G $-equivariant with respect to the diagonal $ G $-actions on $ G \times X $ and $ X \times X. $

\begin{defn}\label{defn_inertia_stack}
	Let $ \cX = [X/G] $ be a quotient stack. The \tb{inertia scheme} of $ X $ under the action of $ G $ is defined by the following pullback diagram:
	\begin{equation*}\label{eq_inertia_scheme}
	\begin{tikzcd}[column sep=4em,row sep=4em]
	I_G X \arrow{r} \arrow{d} & X \arrow{d}{\Delta} \\
	G \times X \arrow{r}{\alpha} & X \times X.
	\end{tikzcd}
	\end{equation*}
	Since $ \Delta $ and $ \alpha $ are $ G $-equivariant, they descend to morphisms $ \Delta: \cX \to \cX \times \cX $ and $ \alpha: [G/G] \times \cX \to \cX \times \cX $. The \tb{inertia stack} of $ \cX $ is defined by the following 2-pullback diagram:
	\begin{equation*}\label{eq_inertia_stack}
	\begin{tikzcd}[column sep=3.6em,row sep=4em]
	I \cX \arrow{r} \arrow{d} & \cX \arrow{d}{\Delta} \\
	\left[G/G\right] \times \cX \arrow{r}{\alpha} & \cX \times \cX.
	\end{tikzcd}
	\end{equation*}
\end{defn}
\begin{rmk}
	The inertia scheme $ I_G X $ can be identified as a subscheme of $ G \times X $ whose $ \CC $-points are
	$$ I_G X (\CC) = \{(g,x) \in G \times X \ | \ gx = x\}. $$
	The inertia stack $ I \cX $ can then be identified as a quotient stack $ [I_G X/G], $ where $ G $ acts on $ I_G X $ by
	$$ g(h,x) = (ghg^{-1}, gx) $$ 
	for $ g, h \in G $ and $ x \in X. $
\end{rmk}

For the rest of this section we will restrict to connected separated quotient DM stacks.
\begin{defn}
	An algebraic stack is connected if its underlying topological space is connected.
\end{defn}
\begin{defn}
	An algebraic stack $ \cX $ is separated if the diagonal morphism $ \Delta: \cX \to \cX \times \cX $ is proper.
\end{defn}
Since properness of morphisms of algebraic spaces is fppf local on the target and preserved under arbitrary base change, we can check the properness of a morphism of algebraic stacks via a fppf representation. For a quotient stack $ \cX = [X/G] $, we have the following 2-pullback diagram:
\begin{equation*}
	\begin{tikzcd}[column sep=3.6em,row sep=4em]
		G \times X \arrow{r} \arrow{d}[swap]{\alpha} & \cX \arrow{d}{\Delta} \\
		X \times X \arrow{r} & \cX \times \cX,
	\end{tikzcd}
\end{equation*}
where the morphism $ X \times X \to \cX \times \cX $ is smooth, and in particular fppf. Therefore, we have an equivalent definition of separatedness of quotient stacks via group actions.
\begin{prop}
	A quotient stack $ \cX = [X/G] $ is separated if and only if the action morphism $ \alpha: G \times X \to X \times X $ is proper.
\end{prop}
\begin{rmk}
	For any separated quotient stack, its stabilizers are finite since they are affine and proper. Therefore, all separated quotient stacks are DM stacks.
\end{rmk}
Let $ \cX = [X/G] $ be a separated quotient DM stack. Note that vector bundles (or coherent sheaves) on $ BG $ are representations of $ G $, and hence $ K^0(BG) = K(BG) = R(G), $ the representation ring of $ G $. Pulling back vector bundles on $ BG $ along the canonical morphism
\begin{equation*}
	q: \cX \to BG
\end{equation*}
and taking their tensor products with coherent sheaves on $ \cX $ makes $ K(\cX) $ a right $ R(G) $-module:
\begin{equation*}
[\cE] \cdot [\rho] = [\cE \otimes q^*\rho]
\end{equation*}
for a coherent sheaf $ \cE $ on $ \cX $ and a representation $ \rho $ of $ G. $
Tensoring by $ \CC, $ we then have a module $ K(\cX)_\CC $ over $ R(G)_\CC. $

Let $ \Supp K(\cX)_\CC $ denote the support of the $ R(G)_\CC $-module $ K(\cX)_\CC $. We will follow the treatment in \cite[Section 4.3]{edidin2013riemann} to decompose the inertia stack $ I\cX $ via $ \Supp K(\cX)_\CC. $ We first state a basic fact which combines two results in \cite[Proposition 2.5]{edidin2005nonabelian} and \cite[Remark 5.1]{edidin2000riemann}.
\begin{lem}
	Let $ G $ be a linear algebraic group. There is a bijection between the set of conjugacy classes of diagonalizable elements in $ G $ and the closed points in $ \Spec R(G)_\CC. $ If $ G $ acts on a separated scheme $ X $ of finite type with finite stabilizers, then the $ R(G)_\CC $-module $ K([X/G])_\CC $ is supported at a finite number of closed points in $ \Spec R(G)_\CC. $
\end{lem}

By the lemma above, $ K(\cX)_\CC $ is supported at a finite number of closed points $ p_0, \dots, p_l $ in $ \Spec R(G)_\CC $ which correspond to a finite number of conjugacy classes $ [g_0] = [1], \dots, [g_l] $ of diagonalizable elements in $ G. $ If $ [g] $ is one of these classes, then the fixed point locus $ X^g $ is non-empty, and hence it must be of finite order. On the contrary, if $ X^g $ is non-empty for an element $ g $, then the corresponding closed point $ \mathfrak{m}_g $ in $ \Spec R(G)_\CC $ is in the support of $ K(\cX)_\CC $. Therefore, we have the following result.
\begin{prop}
	Let $ \cX = [X/G] $ be a separated quotient DM stack. Then the support of the $ R(G)_\CC $-module $ K(\cX)_\CC $ can be identified with a finite set of conjugacy classes of elements of finite order in $ G. $ Furthermore, $ X^g $ is non-empty if and only if $ [g] $ is one of those conjugacy classes.
\end{prop}
Now we assume $ \cX $ is connected. The we have a decomposition
\begin{equation*}
I\cX \simeq \coprod_{i = 0}^{l} \{g_i\} \times \left[X^{g_i} / Z_{g_i}\right],
\end{equation*}
where $ \simeq $ means an equivalence of stacks as categories, $ g_0 = 1 $ is the identity in $ G, $ $ X^{g_i} $ is the fixed point locus in $ X $ under the action of $ g_i, $ and $ Z_{g_i} $ is the centralizer of $ g_i $ in $ G $. Note that each quotient stack $ \{g_i\} \times \left[X^{g_i} / Z_{g_i}\right] $ is not connected in general, but can be further decomposed. We can decompose $ I\cX $ into connected components as follows:
\begin{equation*}
I\cX = \cX_0 \coprod I_t \cX, \quad I_t\cX = \coprod_{i=1}^{l}\coprod_{j=1}^{m_i} \{g_i\} \times [X_{ij}/G_{ij}], % = \coprod_{i,j} \{g_{i}\} \times [X_{ij} /G_{ij}],
\end{equation*}
where $ \cX_0 = \{1\} \times \cX$, $ X_{ij} \subseteq X^{g_i} $ (i.e., $ g_i $ acts on $ X_{ij} $ trivially), and $ G_{ij} \subseteq Z_{g_i}. $ Note that the index $ g_i $ is important, as we will use it to define the orbifold Chern character and the orbifold Todd class later. Re-ordering these components by a single index set $ I $ of size
\begin{equation*}
|I| = 1 + m_1 + \cdots + m_l,
\end{equation*}
we obtain the following.
\begin{prop}\label{prop_decomp_inertia}
	Let $ \cX = [X/G] $ be a connected separated quotient DM stack. Then its inertia stack has a decomposition
	\begin{equation*}
	I\cX \simeq \coprod_{i \in I} \cX_i
	\end{equation*}
	of connected components, where $ i $ runs over a finite index set $ I $ containing $ 0 $, and each $ \cX_i $ is a closed substack of $ \cX $ of the form 
	\begin{equation*}
	\cX_i = \{g_i\} \times [X_i /G_i],
	\end{equation*}
	where $ X_i $ is a closed subscheme of $ X $ invariant under the action of a subgroup $ G_i \subset Z_{g_i} $ and fixed by an element $ g_i \in G_i $ with $ X_0 = X, $ $ G_0 = G $, and $ g_0 = 1. $
\end{prop}
\begin{rmk}
	The finite index set $ I $ in Proposition \ref{prop_decomp_inertia} is called the \tb{inertia index set} of $ \cX. $ The component $ \cX_0 \simeq \cX $ and is called the \tb{distinguished component}, and the other components $ \cX_i $ are called \tb{twisted sectors}.
\end{rmk}

Now we compute two examples of inertia stacks: for the first one, we use the definition of the inertia stack to obtain its decomposition, and for the second one, we use the support of the $ R(G)_\CC $-module $ K(\cX)_\CC $ to obtain the decomposition of $ I\cX. $
\begin{ex}
	Let $ \cX = BG $ be the classifying stack of a finite group $ G $. Then $ I_G X = G $, and the inertia stack of the classifying stack is $ I\cX = [G/G] $, where $ G $ acts on itself by conjugation. Let $ I = \{0, 1, \dots, n\} $ index the conjugacy classes $ [g_0], [g_1], \dots, [g_n] $ of $ G $, where $ g_0 = 1 $ is the identity in $ G $. Then we have a decomposition
	\begin{equation*}
	IBG \simeq \coprod_{i \in I} \{g_i\} \times BZ_{g_i}.
	\end{equation*}
	Here $ \{g_0\} \times BZ_{g_0} = \{1\} \times BG $ is the distinguished component. If $ G $ is abelian, then every element $ g $ forms a conjugacy class and each $ Z_{g} = G, $ so we have $$ IBG = G \times BG. $$
\end{ex}

\begin{ex}
	Let $ X = \CC^2 - \{(0,0)\} $, and let $ G = \CC^* $ act on $ X $ by $ \lambda \cdot (x_0, x_1) = (\lambda^2 x_0, \lambda^3 x_1) $. The quotient stack $ \cX = [X/G] $ is the weighted projective line $ \PP(2,3) $. We have $ K(BC^*)_\CC \cong \CC[x, x^{-1}], $ and
	\begin{equation*}
		K(\PP(2,3))_\CC \cong \frac{\CC[x, x^{-1}]}{(1-x^{-2})(1-x^{-3}))}.
	\end{equation*}
	Therefore, we have
	\begin{equation*}
		\Supp K(\PP(2,3))_\CC = \{g_0, g_1, g_2, g_3\},
	\end{equation*}
	where $ g_0 = 1, g_1 = -1, g_2 = (-1+\sqrt{3}i)/2, $ and $ g_3 = (-1-\sqrt{3}i)/2 $. Now we compute the fixed point loci: $ X^{g_0} = X,\ X^{g_1} = \CC^* \times \{0\}, $ and $\ X^{g_2} = X^{g_3} = \{0\} \times \CC^* $. The centralizer $ Z_{g_i} = \CC^* $ for $ i = 0, 1, 2, 3. $ Therefore, we have a decomposition
	\begin{align*}
	I\PP(2,3) = \coprod_{i=0}^3 \{g_i\} \times [X^{g_i}/Z_{g_i}] \cong \{1\} \times \PP(2,3) \coprod \{-1\} \times B\mu_2 \coprod \left\{\frac{-1 \pm \sqrt{3}i}{2}\right\} \times B\mu_3.
	\end{align*}
\end{ex}
\vskip 2em

\section{Orbifold Chern Character and Orbifold Todd Class}\label{Orbifold_Chern}
In this section we define orbifold Chern characters and orbifold Todd classes of coherent sheaves on a connected separated smooth quotient DM stack.
\vskip 4pt
We first define the Chow ring of a smooth quotient stack following \cite{edidin1998equivariant}.
\begin{defn-prop}
	Let $ \cX = [X/G] $ be a quotient stack. For each integer $ i \geq 0, $ the $ i $-th \tb{Chow group} of $ \cX $ is defined by
	$$ A^i(\cX) = A^i((X \times U)/G), $$ 
	where $ U $ is an open subscheme of codimension $ > i $ of an affine space $ \CC^n $ with a linear $ G $-action such that it restricts to a free $ G $-action on $ U $. The definition of $ A^i(\cX) $
	is independent of the choice of $ n, $ the $ G $-action on $ \CC^n $, and the open subscheme $ U. $ If $ \cX $ is a smooth quotient stack, then there is an intersection product 
	\begin{equation*}
	A^i(\cX) \times A^j(\cX) \to A^{i+j}(\cX)
	\end{equation*}
	for all non-negative integers $ i $ and $ j $ such that $ A^*(\cX) $ becomes a ring graded in $ \ZZ_{\geq 0}. $ We denote by $ A^*(\cX) $ or simply $ A(\cX) $ the \tb{Chow ring} of $ \cX. $ If $ \cX $ is a separated smooth quotient DM stack, then the proper pushforward 
	\begin{equation*}
	A(\cX)_\QQ \to A(X/G)_\QQ,
	\end{equation*}
	where $ X/G $ is the coarse moduli space of $ \cX $,
	is an isomorphism of abelian groups in each degree. This implies that an element of each $ A^i(\cX)_\QQ $ is a linear combination of classes of irreducible closed substacks of $ \cX $ of codimension $ i, $ and $ A^i(\cX)_\QQ $ vanishes for $ i > \dim \cX. $
\end{defn-prop}

\begin{ex}
	The Chow ring of the weighted projective stack $ \PP = \PP(a_0, \dots, a_n) $ is
	\begin{equation*}
	A(\PP) \cong \frac{\ZZ[h]}{(a_0\cdots a_nh^{n+1})},
	\end{equation*}
	where $ h = c_1(\cO_\PP(1)) $ in $ A^1(\PP). $
\end{ex}
\begin{rmk}
	If $ \cX $ is a connected separated smooth quotient DM stack, then its inertia $ I\cX $ is a separated smooth quotient DM stack since the canonical morphism $ I\cX \to \cX $ is finite, and its connected components are smooth because a fixed point locus of a smooth scheme over a field of characteristic zero is smooth.
\end{rmk}
Let $ \cX $ be a connected separated smooth quotient DM stack. There is an orbifold Chern character map
$$ \orbch : K^0(\cX) \to A(I\cX)_{\CC}, $$
which is defined in \cite{toen1999theoremes}.
There is also an orbifold Todd class map
$$ \orbtd : K^0(\cX) \to A(I\cX)_{\CC}^\times, $$
which is given in \cite[Appendix A]{tseng2010orbifold} and \cite[Section 2.4]{iritani2009integral}. Since $ K(\cX) \cong K^0(\cX) $, we have two maps:
$$ \orbch : K(\cX) \to A(I\cX)_{\CC} \quad \text{and} \quad \orbtd : K(\cX) \to A(I\cX)_{\CC}^\times. $$ 

Now we give detailed constructions for these two maps.

\begin{notation}\label{notation_quotient_stack}
	Let $ \cX = [X/G] $ be a connected separated smooth quotient DM stack. Recall the decomposition 
	$$ I\cX \simeq \coprod_{i \in I} \cX_i = \coprod_{i \in I} \{g_i\} \times [X_i/G_i] $$ 
	in Proposition \ref{prop_decomp_inertia}. 
	We have two rings:
	\begin{equation*}
	K^0(I\cX) = \bigoplus_{i \in I} K^0(\cX_i) \quad \text{and} \quad A(I\cX) = \bigoplus_{i \in I} A(\cX_i).
	\end{equation*}
	Let $ \pi: I\cX \to \cX $ denote the natural morphism. For each $ i \in I, $ let $ p_i: X_i \to X $ and $ q_i: \cX_i \to \cX $ denote the inclusions of the closed subscheme $ X_i $ in $ X $ and the closed substack $ \cX_i $ in $ \cX $ respectively, and let $ r_i: \cX_i \to I\cX $ denote the inclusion of $ \cX_i $ in $ I\cX. $ Note that $ \pi \circ r_i = q_i $ for each $ i \in I. $
\end{notation}

We first prove a basic fact that is crucial for the constructions of both the orbifold Chern character and the orbifold Todd class.
\begin{lem}\label{lem_eigenbundle}
	Let $ [X/G] $ be a quotient stack. Let $ g $ be an element in the center of $ G $ such that it has a finite order $ n $ and fixes $ X $. Then every $ G $-equivariant vector bundle $ \cV = (V, \phi) $ on $ X $ admits a $ G $-equivariant decomposition
	\begin{equation}\label{eq_equiv_decomp}
	\cV = \bigoplus_{\lambda \in \eig(\phi(g))} \cV_\lambda,
	\end{equation}
	where $ \eig(\phi(g)) $ consists of eigenvalues of the linear automorphism $ \phi(g) $ on $ V $ which are $ n $th roots of unity, and each $ \cV_\lambda = (V_\lambda, \phi_\lambda) $ is a $ G $-equivariant vector bundle on $ X. $
\end{lem}
\begin{proof}
	Consider a $ G $-equivariant vector bundle $ (V, \phi) $ on $ X. $ This corresponds to an action
	\begin{equation*}
	\phi: G \to \Aut(V)
	\end{equation*}
	of $ G $ on $ V. $ By assumption, $ g^n = 1 $, so the linear automorphism $ \phi(g) $ on $ V $ satisfies $ (\phi(g))^n = \id_V, $ and hence each eigenvalue of $ \phi(g) $ is an $ n $th root. Since $ g $ acts on $ X $ trivially, there is an eigenbundle decomposition
	$$ V = \bigoplus_{\lambda \in \eig(\phi(g))} V_\lambda $$
	such that for each $ v \in V_\lambda, $
	\begin{equation*}
	\phi(g) (v) = \lambda v.
	\end{equation*}
	Take any $ h \in G. $ Since $ g $ is in the center of $ G, $ we have $ gh = hg. $ For any $ v \in V_\lambda, $ we have
	\begin{equation*}
	\phi(g) \left(\phi(h)v\right) = \phi(h) \left(\phi(g)v\right) = \phi(h) \left(\lambda v\right) = \lambda \left(\phi(h)v\right),
	\end{equation*}
	which implies each $ V_\lambda $ is an invariant subbundle of $ V $ under the action of $ G $, i.e., $ \phi $ restricts to an action
	\begin{equation*}
	\phi_\lambda: G \to \Aut(V_\lambda)
	\end{equation*}
	of $ G $ on each $ V_\lambda $. Therefore, we have the $ G $-equivariant decomposition of $ (V, \phi) $ claimed in the lemma.
\end{proof}

\begin{defn}
	Let $ \cX = [X/G] $ be a quotient stack. Let $ g $ be an element in the center of $ G $ such that it has a finite order and fixes $ X $.
	The $ \boldsymbol{g} $\tb{-twisting morphism}
	\begin{equation*}
	\rho_g: K^0(\cX) \to K^0(\cX)_\CC
	\end{equation*}
	is defined by
	\begin{equation*}
	\rho_g ([\cV]) = \sum_{\lambda \in \eig(\phi(g))} \lambda [\cV_\lambda] 
	\end{equation*}
	for the decomposition (\ref{eq_equiv_decomp}) of a vector bundle $ \cV $ on $ \cX $ and extended linearly. 
\end{defn}
The following lemma is immediate.
\begin{lem}\label{lem_twisting_morphism}
	The $ g $-twisting morphism $ \rho_g: K^0(\cX) \to K^0(\cX)_\CC $ is a ring homomorphism.
\end{lem}
\begin{defn}\label{defn_twisted_Chern}
	Let $ \cX = [X/G] $ be a separated smooth quotient stack. Let $ g $ be an element in the center of $ G $ such that it has a finite order and fixes $ X $. We define the $ \boldsymbol{g} $\tb{-twisted Chern character} map by
	\begin{equation*}
	\ch^{\rho_g}: K^0(\cX) \xrightarrow{\rho_g} K^0(\cX)_\CC \xrightarrow{\ch} A(\cX)_\CC,
	\end{equation*}
	where $ \ch: K^0(\cX)_\CC \to A(\cX)_\CC $ is the Chern character map on the stack $ \cX $ tensored by $ \CC. $
\end{defn}
Lemma \ref{lem_twisting_morphism} implies the following.
\begin{prop}\label{lem_twisted_Chern}
	The $ g $-twisted Chern character map $ \ch^{\rho_g}: K^0(\cX) \to A(\cX)_\CC $ is a ring homomorphism.
\end{prop}
\begin{ex}\label{ex_twisted_Chern_on_BG}
	Consider the classifying stack $ BG $ for a finite group $ G. $ Then $ K^0(BG) = K(BG) $ is the representation ring $ R(G). $ Take an element $ g $ in the center of $ G. $ Let's compute the $ g $-twisted Chern character
	\begin{equation*}
	\ch^{\rho_{g}}: R(G) \to A(BG)_\CC \cong \CC.
	\end{equation*}
	A representation $ \phi: G \to \GL(V) $ of $ G $ on a vector space $ V $ decomposes as
	\begin{equation*}
	\phi = \bigoplus_{\lambda \in \eig(\phi(g))} \phi_\lambda
	\end{equation*}
	with subrepresentations $ \phi_\lambda: G \to \GL(V_\lambda), $ so we have a $ g $-twisting morphism
	\begin{equation*}
	\rho_{g}: R(G) \to R(G)_\CC, \quad \phi \mapsto \sum_{\lambda \in \eig(\phi(g))} \lambda \phi_\lambda.
	\end{equation*}
	The Chern character map $ \ch: R(G)_\CC \to \CC $ is a linear extension of the rank map, i.e., for a virtual representation $ \sum_i a_i \vphi_i $ of $ G $, where $ a_i \in \CC $ and representation $ \vphi_i: G \to \GL(V_i) $, we have
	\begin{equation*}
	\ch\left(\sum_i a_i \vphi_i\right) = \sum_i a_i \dim V_i.
	\end{equation*}
	Therefore, the $ g $-twisted Chern character map is given by 
	\begin{equation*}
	\ch^{\rho_{g}}(\phi) = \ch\left(\rho_g(\phi)\right) = \sum_{\lambda \in \eig(\phi(g))} \lambda \dim V_\lambda = \tr(\phi(g)) = \chi_\phi(g),
	\end{equation*}
	where $ \chi_\phi $ denotes the character of the representation $ \phi $ of $ G. $ The fact that the twisted Chern character map is a ring map implies that for any pair of representations $ \phi $ and $ \psi $ of $ G, $ we have
	\begin{equation*}
	\chi_{\phi \oplus \psi}(g) = \chi_{\phi}(g) + \chi_{\psi}(g) \quad \text{and} \quad \chi_{\phi \otimes \psi}(g) = \chi_{\phi}(g) \chi_{\psi}(g),
	\end{equation*}
	which is a basic fact in character theory.
\end{ex}
\begin{defn-prop}
	Let $ \cX = [X/G] $ be a connected separated smooth quotient DM stack. We have a finite inertia index set $ I $ and a decomposition 
	$$ I\cX \simeq \coprod_{i \in I} \cX_i = \coprod_{i \in I} \{g_i\} \times [X_i/G_i]. $$ 
	Each $ g_i $ is in the center of $ G_i $ (because $ G_i $ is contained in the centralizer $ Z_{g_i} $ of $ g_i $), it fixes $ X_i $ and has a finite order, so it induces a $ g_i $-twisting morphism
	\begin{equation*}
		\rho_{g_i}: K^0(\cX_i) \to K^0(\cX_i)_\CC,
	\end{equation*}
	and hence a $ g_i $-twisted Chern character map
	\begin{equation*}
		\ch^{\rho_{g_i}}: K^0(\cX_i) \to A(\cX_i)_\CC.
	\end{equation*}
	Let $ \rho: K^0(I\cX) \to K^0(I\cX)_\CC $ denote the \tb{twisting morphism} on $ K^0(I\cX) $, i.e.,
	\begin{equation*}
		\rho = \bigoplus_{i \in I} \rho_{g_i}.
	\end{equation*}
\end{defn-prop}

Now we have all the ingredients to define the orbifold Chern character map.
\begin{defn}\label{defn_orbifold_Chern}
	Let $ \cX = [X/G] $ be a connected separated smooth quotient DM stack.
	The \tb{inertia Chern character} map is defined by the composition
	$$ \ch^\rho: K^0(I\cX) \xrightarrow{\rho} K^0(I\cX)_\CC \xrightarrow{\ch} A(I\cX)_\CC. $$
	The \tb{orbifold Chern character} map on $ K(\cX) $ is defined by the composition
	\begin{equation*}
	\orbch : K(\cX) \xrightarrow{\beta} K^0(\cX) \xrightarrow{\pi^*} K^0(I\cX) \xrightarrow{\ch^\rho} A(I\cX)_\CC.
	\end{equation*}
\end{defn}
\begin{prop}\label{prop_orbch_ring_map}
	The orbifold Chern character map $ \orbch : K(\cX) \to A(I\cX)_{\CC} $ is a ring homomorphism, i.e.,
	\begin{equation*}
	\quad \orbch(x+y) = \orbch(x) + \orbch(y) \quad \text{and} \quad \orbch(xy) = \orbch(x) \orbch(y)
	\end{equation*}
	for all $ x $ and $ y $ in $ K(\cX). $
\end{prop}
\begin{proof}
	Lemma \ref{lem_twisted_Chern} implies the inertia Chern character map $ \ch^\rho: K^0(I\cX) \to  A(I\cX)_\CC $ is a ring homomorphism, where both $ K^0(I\cX) $ and $ A(I\cX)_\CC $ carry ring structures from the direct sum of the rings $ K^0(\cX_i) $ and $ A(\cX_i)_\CC $ over $ I. $ The orbifold Chern character $ \orbch : K(\cX) \to A(I\cX)_{\CC} $ is a ring homomorphism since it's a composition of three ring homomorphisms $ \beta, \pi^* $ and $ \ch^\rho $.
\end{proof}
\begin{rmk}[Explicit formulas for the orbifold Chern character]
	The orbifold Chern character of a vector bundle $ \cV $ on $ \cX $ is given by
	\begin{equation*}\label{eq_orbch_concrete}
	\orbch (\cV) = \left(\ch(\cV),\ \bigoplus_{0 \neq i \in I} \ch^{\rho_{g_i}}(\cV_i)\right) = \left(\ch(\cV), \  \bigoplus_{0 \neq i \in I} \sum_{\lambda \in \eig(\phi_i(g_i))} \lambda \ch (\cV_{i,\lambda})\right), \vspace{5pt}
	\end{equation*}
	where each $ \cV_i $ is the restriction of $ \cV $ from $ \cX $ to the substack $ \cX_i, $ with an eigenbundle decomposition $ \cV_i = \oplus_{\lambda \in \eig(\phi_i(g_i))} \cV_{i,\lambda} $ for $ 0 \neq i \in I $.
\end{rmk}

\begin{rmk}[Care is required when working with sheaves.]
	Consider two coherent sheaves $ \cE, \cF $ on $ \cX. $ Recall from Proposition \ref{lem_multiplication_sheaves} that their product is given by
	\begin{equation*}
	[\cE] [\cF] = \sum_i (-1)^i [\cE_i \otimes \cF]
	\end{equation*}
	for any finite resolution $ \cE_{\boldsymbol{\cdot}} \to \cE \to 0. $
	The identity 
	$$ \orbch([\cE] [\cF]) = \orbch([\cE]) \orbch([\cF]) $$ 
	implies
	\begin{equation*}
	\orbch(\cE) \orbch(\cF) = \sum_i (-1)^i \orbch(\cE_i \otimes \cF),
	\end{equation*}
	but
	\begin{equation*}
	\orbch(\cE) \orbch(\cF) \neq \orbch(\cE \otimes \cF)
	\end{equation*}
	in general, since $ \cE_{\boldsymbol{\cdot}} \otimes \cF \to \cE \otimes \cF \to 0 $ may fail to be exact. However, if either $ \cE $ or $ \cF $ is locally free, then we have
	\begin{equation*}
	\orbch(\cE) \orbch(\cF) = \orbch(\cE \otimes \cF).
	\end{equation*}
\end{rmk}
Recall that for a separated smooth scheme $ X $ of finite type, the Chern character map $ \ch: K(X)_\QQ \xrightarrow{\sim} A(X)_\QQ $
is a $ \QQ $-algebra isomorphism after tensored with $ \QQ. $ A similar result holds for quotient DM stacks.
\begin{prop}\label{prop_orbch_isom}
	Let $ \cX $ be a connected separated smooth quotient DM stack. The orbifold Chern character map
	$$ \orbch : K(\cX)_\CC \xrightarrow{\sim} A(I\cX)_{\CC} $$
	is a $ \CC $-algebra isomorphism after tensored with $ \CC. $
\end{prop}
\begin{proof}
	Note that $ \cX $ is a separated DM stack of finite type. By the discussion on page 9 in \cite{toen2000motives}, after tensored with $ \QQ, $ the orbifold Chern character $ \orbch: K(\cX)_\QQ \to A_\chi(\cX) $ is a $ \QQ $-algebra isomorphism onto $ A_\chi(\cX) $, the rational Chow ring with coefficients in the characters of $ \cX $. Since $ A_\chi(\cX)_\CC \cong A(I\cX)_{\CC}, $ the map $ \orbch $ becomes a $ \CC $-algebra isomorphism after tensored with $ \CC. $
\end{proof}
\begin{rmk}
	The ordinary Chern character $ \ch: K^0(\cX)_\CC \to A(\cX)_\CC  $ is a surjective $ \CC $-algebra homomorphism, but not an isomorphism in general, as the following example shows.
\end{rmk}
\begin{ex}[The orbifold Chern character for $ BG $ is the character of representations of $ G $.]\label{ex_orbch_BG}
	Let $ \cX = BG $ for a finite group $ G $ with $ n $ conjugacy classes 
	$$ [g_0], [g_1], \dots, [g_{n-1}]. $$
	Recall that 
	$$ IBG = \coprod_{i=0}^{n-1} \{g_i\} \times BZ_{g_i}. $$ 
	Take a vector bundle $ \cV = (V, \phi) $ on $ BG, $ i.e., a representation $\phi: G \to \GL(V) $
	of $ G $. 
	On each component $ \{g_i\} \times BZ_{g_i} $, the representation $ \phi $ of $ G $ restricts to a representation
	\begin{equation*}
	\phi_i: Z_{g_i} \to \GL(V)
	\end{equation*}
	of $ Z_{g_i} $.
	Since each $ g_i $ is in the center of $ Z_{g_i}, $ by the computations in Example \ref{ex_twisted_Chern_on_BG}, we know that each $ g_i $-twisted Chern character map is
	\begin{equation*}
	\ch^{\rho_{g_i}}: R(Z_{g_i}) \to \CC, \quad \vphi \mapsto \chi_\vphi(g_i).
	\end{equation*}
	Therefore, the orbifold Chern character becomes
	\begin{alignat*}{3}
	\orbch: & R(G) & \ \longrightarrow & \ \bigoplus\limits_{i=0}^{n-1} R(Z_{g_i}) & \ \longrightarrow & \ \CC^n \\
	& \quad \phi & \ \longmapsto & \ (\phi_0, \dots, \phi_{n-1}) & \ \longmapsto & \ \left( \chi_\phi(g_0), \cdots, \chi_\phi(g_{n-1}) \right),
	\end{alignat*}
	where each $ \chi_\phi(g_i) = \chi_{\phi_i}(g_i) $ since $ \phi_i $ is a restriction of $ \phi $.
	Observe that the zeroth component 
	$$ \chi_\phi(g_0) = \chi_\phi(1) = \deg \phi = \dim V $$ 
	is the ordinary Chern character of $ \phi. $
	The space $ \CC^n $ can be identified with the space $ C(G) $ of class functions on $ G. $
	Therefore, the orbifold Chern character coincides with the character map
	\begin{align*}
	\orbch: R(G) & \rightarrow C(G) \cong \CC^n \\
	\phi & \mapsto \chi_\phi.
	\end{align*}
	This is an injective ring map, where irreducible representations map to irreducible characters which form a basis of $ C(G) $. Tensoring $ R(G) $ with $ \CC $, we obtain a $ \CC $-algebra isomorphism
	$$ \orbch: R(G)_\CC \xrightarrow{\sim} C(G). $$
	Compare with the ordinary Chern character map (i.e., the rank map)
	\begin{align*}
	\ch: R(G) \rightarrow \ZZ, \quad (V, \phi) \mapsto \dim V,
	\end{align*}
	which is surjective but not injective in general.
\end{ex}
\begin{ex}[The orbifold Chern character map on $ B\mu_n $ is the inverse discrete Fourier transform]\label{ex_discrete_FT}
	Take $ G = \mu_n $ in Example \ref{ex_orbch_BG}. Then its representation ring $ R(\mu_n) \cong \ZZ[x]/(x^n-1) $,
	where $ x $ is the character
	$$ \mu_n \to \CC^*, \quad \omega = e^{2\pi i/n} \mapsto \omega $$
	of $ \mu_n. $ The orbifold Chern character for $ \mu_n $ is the ring map
	\begin{align*}
	\orbch: \frac{\ZZ[x]}{(x^n-1)} \to \CC^n, \quad x \mapsto (1, \omega, \dots, \omega^{n-1}).
	\end{align*}
	Tensoring $ R(\mu) $ with $ \CC, $ we obtain a $ \CC $-algebra isomorphism 
	$$ \orbch: \frac{\CC[x]}{(x^n-1)} \xrightarrow{\sim} \CC^n. $$
	Here the inverse map $ \orbch^{-1}: \CC^n \to \CC[x]/(x^n-1) $ can be identified with the discrete Fourier transform as follows. For each $ k = 0, \dots, n-1, $ put
	\begin{equation*}
	e_k = \orbch(x^k) = (1, \omega^k, \dots, \omega^{(n-1)k}).
	\end{equation*}
	Consider the weighted inner product on $ \CC^n $ defined by
	\begin{equation*}
	\inprod{a,b}_W = \frac{1}{n} \sum_{i=0}^{n-1} \cj{a}_i b_i
	\end{equation*}
	for all $ a = (a_0, \dots, a_{n-1}) $ and $ b = (b_0, \dots, b_{n-1}) $ in $ \CC^n. $ Then the elements $ e_0, \dots, e_{n-1} $ form an orthonormal basis, i.e., for all $ i,j = 0, \dots, n-1, $ we have
	\begin{equation*}
	\inprod{e_i,e_j}_W = \delta_{ij},
	\end{equation*}
	where $ \delta_{ij} $ is the Kronecker delta. Take any function $ f: \ZZ/n\ZZ \to \CC. $ Identifying $ f $ with an element $ (f(0), \dots, f(n-1)) $ in $ \CC^n, $ we then have the orthogonal decomposition
	\begin{equation*}
	f = \sum_{k=0}^{n-1} \inprod{e_k, f}_W e_k, 
	\end{equation*}
	where
	\begin{equation*}
	\inprod{e_k, f}_W = \frac{1}{n} \sum_{j=0}^{n-1} \cj{\omega}^{jk} f(j) = \frac{1}{n} \sum_{j=0}^{n-1} e^{-2\pi ijk/n} f(j)
	\end{equation*}
	% in $ \CC^n $ via the orthonormal basis $ \{e_0, \dots, e_{n-1}\}. $ 
	for $ k = 0, \dots, n-1. $ This defines a map
	\begin{equation*}
	\psi: \CC^n \to \frac{\CC[x]}{(x^n-1)}, \quad f \mapsto \hat{f} = \sum_{k=0}^{n-1} \inprod{e_k, f}_W x^k,
	\end{equation*}
	which can be identified with the \tb{discrete Fourier transform} (up to a scaling depending on  the convention). The map $ \psi $ is indeed the inverse of $ \orbch $ since $ \psi(e_k) = x^k $ for all $ k = 0, \dots, n-1. $ Therefore, the orbifold Chern character for $ \mu_n $ is nothing but the inverse discrete Fourier transform, which maps an element
	\begin{equation*}
	f = f(0) + f(1) x + f(2) x^2 + \dots + f(n-1) x^{n-1}
	\end{equation*}
	in $ \CC[x]/(x^n-1) $ to an element $ \check{f} = (\check{f}(0), \dots, \check{f}(n-1)) $ in $ \CC^{n} $, where
	\begin{equation*}
	\check{f}(k) = \sum_{j=0}^{n-1} \omega^{jk} f(j) = \sum_{j=0}^{n-1} e^{2\pi ijk/n} f(j)
	\end{equation*}
	for $ k = 0, \dots, n-1. $
\end{ex}

Next we will define the orbifold Todd class map.
\begin{defn}
	Let $ \cX = [X/G] $ be a quotient stack. Let $ g $ be an element in the center of $ G $ such that it has a finite order and fixes $ X $. For the decomposition (\ref{eq_equiv_decomp}) of a vector bundle $ \cV = (V, \phi) $ on $ \cX, $ we write
	\begin{equation*}
	\cV^\fix = \cV_1 \quad \text{and} \quad \cV^\mov = \bigoplus_{1 \neq \lambda \in \eig(\phi(g))} \cV_\lambda.
	\end{equation*}
	We define two projections
	\begin{equation*}
	P_f: K^0(\cX) \to K^0(\cX) \quad \text{and} \quad P_m: K^0(\cX) \to K^0(\cX)
	\end{equation*}
	by setting
	\begin{equation*}
	P_f(x) = \cV^\fix - \cW^\fix \quad \text{and} \quad P_m(x) = \cV^\mov - \cW^\mov
	\end{equation*}
	for any element $ x = [\cV] - [\cW] $ in $ K^0(\cX) $, where $ \cV $ and $ \cW $ are vector bundles on $ \cX. $
\end{defn}
\begin{rmk}
	It's easy to check both $ P_f $ and $ P_m $ are indeed projections, i.e., they are linear and idempotent.
\end{rmk}
\begin{defn}
	Let $ \cX = [X/G] $ be a connected separated smooth quotient DM stack. Let $ g $ be an element in the center of $ G $ such that it has a finite order and fixes $ X $. The $ \boldsymbol{g} $\tb{-twisted Euler class} map is defined by
	\begin{equation*}
	e^{\rho_g}: P_m(K^0(\cX)) \xrightarrow{e^K} K^0(\cX) \xrightarrow{\ch^{\rho_g}} A(\cX)_\CC,
	\end{equation*}
	where $ P_m(K^0(\cX)) $ is the image of the projection $ P_m. $
\end{defn}
\begin{prop}
	The $ g $-twisted Euler class $ e^{\rho_g} $ is multiplicative and maps into the units of $ A(\cX)_\CC, $ i.e., we have a group homomorphism
	\begin{equation*}
	e^{\rho_g}: \left(P_m(K^0(\cX)), +\right) \to \left(A(\cX)_\CC^\times, \times\right).
	\end{equation*}
\end{prop}
\begin{proof}
	It suffices to consider an element $ x = P_m(\cV) $ in $ K^0(\cX) $ for a vector bundle $ \cV = (V, \phi) $ on $ \cX. $
	By the splitting principle, we can assume that each eigenbundle $ \cV_{\lambda} = (V_\lambda, \phi_\lambda) $ decomposes into line bundles 
	\begin{equation*}
	\cV_{\lambda} = \bigoplus_{j=1}^{\rk(\cV_{\lambda})} \cL_{\lambda,j}.
	\end{equation*}
	Let $ x_{\lambda,j} = [\cL_{\lambda,j}] $ in $ K^0(\cX) $. Then we have
	\begin{equation*}
	x = \sum_{1 \neq \lambda \in \eig(\phi(g))} \sum_{j=1}^{\rk(\cV_{\lambda})} x_{\lambda,j}
	\end{equation*}
	in $ K^0(\cX) $. Since $ e^K $ is multiplicative, we have
	\begin{align*}
	e^{\rho_g}(x) & = \ch^{\rho_g} \left(\prod_{\lambda \neq 1,j} e^K(x_{\lambda,j})\right) = \left(\ch \circ \rho_g\right) \prod_{\lambda \neq 1,j} \left(1-x_{\lambda,j}^\vee\right) \\
	& = \ch \left(\prod_{\lambda \neq 1,j} \left( 1 - \lambda^{-1} x_{\lambda,j}^\vee \right)\right) = \prod_{\lambda \neq 1,j} \left( 1 - \lambda^{-1} e^{-h_{\lambda,j}} \right)
	\end{align*}
	in $ A(\cX)_\CC, $ where $ h_{\lambda,j} = c_1(x_{\lambda,j}). $ Since the degree zero component in $ e^{\rho_g}(x) $ is nonzero, it is invertible in $ A(\cX)_\CC. $
\end{proof}
\begin{defn}\label{defn_twisted_Todd}
	Let $ \cX = [X/G] $ be a connected separated smooth quotient DM stack. Let $ g $ be an element in the center of $ G $ such that it has a finite order and fixes $ X $. We define the $ \boldsymbol{g} $\tb{-twisted Todd class} map 
	$$ \td^{\rho_g}: K^0(\cX) \to A(\cX)_\CC^\times $$
	by setting
	\begin{equation*}
	\td^{\rho_g}(x) = \frac{\td(P_f(x))}{e^{\rho_g}(P_m(x))}
	\end{equation*}
	for any $ x $ in $ K^0(\cX). $
\end{defn}
\begin{rmk}
	Applying Definition \ref{defn_twisted_Todd} to each substack $ \cX_i = \{g_i\} \times [X_i/G_i] $ of $ \cX $, we then have a $ g_i $-twisted Todd class map
	\begin{equation*}
		\td^{\rho_{g_i}}: K^0(\cX_i) \to A(\cX_i)_\CC^\times.
	\end{equation*}
\end{rmk}

Now we can define the orbifold Todd class map.
\begin{defn}
	Let $ \cX = [X/G] $ be a connected separated smooth quotient stack in Notation \ref{notation_quotient_stack}. Define projections $ IP_f $ and $ IP_m $ on $ K^0(I\cX) $ by the direct sums
	\begin{equation*}
	IP_f = \bigoplus_{i \in I} P_f \quad \text{and} \quad IP_m = \bigoplus_{i \in I} P_m.
	\end{equation*}
	Define the \tb{inertia Euler class} map by the composition
	\begin{equation*}
	e^\rho: IP_m(K^0(I\cX)) \xrightarrow{e^K} K^0(I\cX) \xrightarrow{\ch^\rho} A(I\cX)_\CC^\times.
	\end{equation*}
	The \tb{inertia Todd class} map
	$$ \td^\rho: K^0(I\cX) \to A(I\cX)_{\CC}^\times $$
	is defined by
	\begin{equation*}
	\td^\rho(x) = \frac{\td(IP_f(x))}{e^\rho(IP_m(x))}
	\end{equation*}
	for any $ x $ in $ K^0(I\cX) $. The \textbf{orbifold Todd class} map on $ K(\cX) $ is defined by the composition
	\begin{equation*}
	\orbtd : K(\cX) \xrightarrow{\beta} K^0(\cX) \xrightarrow{\pi^*} K^0(I\cX) \xrightarrow{\td^\rho} A(I\cX)_\CC^\times.
	\end{equation*}
\end{defn}
The following result is immediate.
\begin{prop}
	The orbifold Todd class map $ \orbtd : K(\cX) \to A(I\cX)_\CC^\times $ is multiplicative, i.e.,
	\begin{equation*}
	\orbtd(x+y) = \orbtd(x) \orbtd(y)
	\end{equation*}
	for all $ x $ and $ y $ in $ K(\cX). $
\end{prop}
\begin{rmk}[Explicit formulas for the orbifold Todd class]
	The orbifold Todd class of a vector bundle $ \cV = (V, \phi) $ on $ \cX $ is given by
	\begin{align*}
	\orbtd (\cV) = \left(\td(\cV),\ \bigoplus_{0 \neq i \in I} \frac{\td(\cV_i^\fix)}{e^{\rho_{g_i}}(\cV_i^\mov)}\right) = \left(\td(\cV), \ \bigoplus_{0 \neq i \in I}  \frac{\td(\cV_i^\fix)}{\prod_{\lambda \neq 1, j } (1 - \lambda^{-1} e^{-h_{i,\lambda,j}})} \right), \vspace{5pt}
	\end{align*}
	where $ \lambda \in \eig(\phi_i(g_i)) $, and $ h_{i,\lambda,j} $'s are the Chern roots of each eigenbundle $ \cV_{i,\lambda} $ on $ \cX_i $ for $ 0 \neq i \in I $ and $ 1 \leq j \leq \rk(\cV_{i,\lambda}) $.
\end{rmk}
We are interested in the orbifold Todd class of the $ K $-theory class corresponding to the tangent complex of $ \cX. $
\begin{defn}[The tangent complex of a quotient stack]
	Let $ \cX = [X/G] $ be a smooth quotient stack. The tangent bundle $ TX $ has a canonical $ G $-equivariant structure $ \tau $ defined as follows. Take any $ g \in G. $ The map $ x \mapsto gx $ gives a smooth automorphism on $ X $
	\begin{equation*}
	g: X \to X
	\end{equation*}
	which induces a linear isomorphism on the tangent spaces
	\begin{equation*}
	dg_x: T_x X \to T_{gx} X
	\end{equation*}
	at all points $ x \in X. $ Therefore, we have a $ G $-action
	\begin{equation*}
	\tau: G \to \Aut(TX), \quad g \mapsto (\tau(g): (x,v) \mapsto (gx,dg_x(v)).		
	\end{equation*}
	The $ G $-equivariant vector bundle $ (TX, \tau) $ corresponds to a vector bundle on $ \cX, $ which we may simply denote by $ TX. $ Take any $ x \in X. $ The $ G $-action on $ X $ also gives a smooth morphism
	\begin{equation*}
	\sigma_x: G \to X, \quad g \mapsto gx
	\end{equation*}
	which induces a linear map 
	\begin{equation}\label{eq_lieG_Tx}
	L_x : \mathfrak{g} \to T_x X, \quad A \mapsto (d\sigma_x)_1(A),
	\end{equation}
	where $ \mathfrak{g} $ is the Lie algebra of $ G. $ There is a \tb{fundamental vector field} $ s_A $ associated to every $ A \in \mathfrak{g}, $ i.e., a section $ s_A: X \to TX $ with values 
	\begin{equation*}
	s_A(x) = L_x (A)
	\end{equation*}
	in $ T_xX $ for every $ x \in X. $ The \tb{tangent complex} of $ \cX $ is defined by a 2-term complex
	\begin{equation*}
	T\cX = (\cO_X \otimes \mathfrak{g},\phi_\Ad) \xrightarrow{L} (TX, \tau)
	\end{equation*}
	concentrated in degree $ -1 $ and $ 0, $ where $ \Ad $ is the adjoint representation 
	\begin{equation*}
	\Ad: G \to \GL(\fg), \quad g \mapsto (\Ad(g): A \mapsto gAg^{-1})
	\end{equation*}
	of $ G $ and the morphism $ L $ is induced by the linear maps $ L_x $ in (\ref{eq_lieG_Tx}), i.e., $ L $ corresponds to the morphism
	\begin{equation*}
	X \times \fg \to \Tot(TX), \quad (x,A) \mapsto (x,s_A(x))
	\end{equation*}
	between the total spaces of $ \cO_X \otimes \mathfrak{g} $ and $ TX $ which is $ G $-equivariant, as can be easily checked. To simplify notations, we may denote $ T\cX $ by
	\begin{equation*}
	T\cX = \cO_X \otimes \mathfrak{g} \xrightarrow{L} TX.
	\end{equation*}
	The tangent space of $ \cX $ at a point $ x \in \cX $ is given by
	\begin{equation*}
	T_x\cX = T_xX / L_x(\fg).
	\end{equation*}
	Note that the image $ L_x(\fg) \subset T_xX $ can be identified with the tangent space of the orbit $ Gx $ at $ x. $ When $ G $ is a finite group, then $ T\cX = TX $ is simply the $ G $-equivariant tangent bundle on $ X, $ and the tangent space $ T_x \cX = T_x X $ at all $ x \in X. $ 
\end{defn}

Now let's compute the orbifold Todd class of the $ K $-theory class 
\begin{equation*}
[T\cX] = [TX] - [\cO_X \otimes \mathfrak{g}] = [TX, \tau] - [\cO_X \otimes \mathfrak{g}, \phi_\Ad]
\end{equation*}
of the tangent complex $ T\cX $ in $ K(\cX). $ 
By the definition of $ \orbtd $, we have
\begin{equation*}
\orbtd\left(T\cX\right) = \frac{\td(IP_f[\pi^* T\cX])}{e^\rho(IP_m[\pi^* T\cX])},
\end{equation*}
where $ IP_f[\pi^* T\cX] $ and $ IP_m[\pi^* T\cX] $ are computed in the following lemma.
\begin{lem}
	For a connected separated smooth quotient DM stack $ \cX, $ there is a splitting short exact sequence
	\begin{equation*}
	0 \to TI\cX \to \pi^*T\cX \to NI\cX \to 0
	\end{equation*}
	of 2-term complexes on the inertia stack $ I\cX $ such that
	\begin{equation*}
	[TI\cX] = IP_f[\pi^* T\cX] \quad \text{and} \quad [NI\cX] = IP_m[\pi^* T\cX]
	\end{equation*}
	in $ K^0(I\cX). $
\end{lem}
\begin{proof}
	Let $ I $ be the inertia index set of $ \cX $ in Notation \ref{notation_quotient_stack}. Take $ i \in I. $ Pulling back $ T\cX $ along $ q_i: \cX_i \to \cX $ gives a complex
	\begin{equation*}
	(T\cX)_i = q_i^*T\cX = (\cO_{X_i} \otimes \mathfrak{g},(\phi_\Ad)_i) \to ((TX)_i, \tau_i),
	\end{equation*}
	where $ (\phi_\Ad)_i $ and $ \tau_i $ are $ G_i $-equivariant structures on $  \cO_{X_i} \otimes \mathfrak{g} $ and $ (TX)_i = p_i^*TX $ induced from $ G$-equivariant structures $ \phi_\Ad $ and $ \tau $ respectively.
	We want to find the fixed and moved sub-complexes 
	\begin{equation*}
	(T\cX)_i^\fix \quad \text{and} \quad (T\cX)_i^\mov
	\end{equation*}
	of the tangent complex $ T\cX $ under the action of $ g_i $ (via linear operators $ \phi_\Ad(g_i) $ and $ \tau(g_i) $ on the two terms of $ T\cX $). To simplify notations, we write
	\begin{equation*}
	(T\cX)_i = \cO_{X_i} \otimes \mathfrak{g} \to (TX)_i,
	\end{equation*}
	where $ G_i $-equivariant structures are understood. 
	The inclusion $ G_i \into G $ induces a splitting short exact sequence
	\begin{equation*}
	0 \to \fg_i \to \fg \to \fg/\fg_i \to 0
	\end{equation*}
	of representations of $ G_i $, where $ \fg_i $ carries the adjoint representation $$ \Ad_i: G_i \to \GL(\fg_i) $$
	of $ G_i $, and $ \fg $ carries the restriction $ \Ad|_{G_i} $ of the adjoint representation of $ G. $
	Since $ G_i $ is a subgroup of the centralizer $ Z_{g_i}, $ we have
	\begin{equation*}
	\Ad_i(g_i) = \id_{\fg_i}.
	\end{equation*}
	Therefore, we have a splitting short exact sequence
	\begin{equation*}
	0 \to \cO_{X_i} \otimes \fg_i \to \cO_{X_i} \otimes \fg \to \cO_{X_i} \otimes \fg/\fg_i \to 0
	\end{equation*}
	of $ G_i $-equivariant vector bundles on $ X_i, $ where
	\begin{equation*}
	\cO_{X_i} \otimes \fg_i = \left(\cO_{X_i} \otimes \fg\right)^\fix \quad \text{and} \quad \cO_{X_i} \otimes \fg/\fg_i = \left(\cO_{X_i} \otimes \fg\right)^\mov.
	\end{equation*}
	The inclusion $ p_i: X_i \into X $ gives a splitting short exact sequence
	\begin{equation*}
	0 \to TX_i \to p_i^* TX \to NX_i \to 0
	\end{equation*}
	of $ G_i $-equivariant vector bundles on $ X_i, $ where $ TX_i $ and $ NX_i $ are the tangent bundle and normal bundle of $ X_i $ with their canonical $ G_i $-equivariant structures. Since $ g_i $ acts trivially on $ X_i, $ we can identify
	\begin{equation*}
	TX_i = (TX)_i^\fix \quad \text{and} \quad NX_i = (TX)_i^\mov.
	\end{equation*}
	Therefore, we obtain a splitting short exact sequence 
	\begin{equation*}
	0 \to T\cX_i \to (T\cX)_i \to N\cX_i \to 0
	\end{equation*}
	of 2-term complexes on $ \cX_i $, where 
	\begin{equation*}
	T\cX_i = \cO_{X_i} \otimes \fg_i \to TX_i \quad \text{and} \quad N\cX_i = \cO_{X_i} \otimes \fg/\fg_i \to NX_i
	\end{equation*}
	are identified with $ (T\cX)_i^\fix $ and $ (T\cX)_i^\mov $ respectively.
	Define two complexes
	\begin{equation*}
	TI\cX = \coprod_{i \in I} T\cX_i \quad \text{and} \quad NI\cX = \coprod_{i \in I} N\cX_i
	\end{equation*}
	on $ I\cX. $ Since $ \pi^*T\cX = \coprod_{i \in I} \left(T\cX\right)_i, $ we have
	\begin{equation*}
	[TI\cX] = IP_f[\pi^* T\cX] \quad \text{and} \quad [NI\cX] = IP_m[\pi^* T\cX]
	\end{equation*}
	in $ K^0(I\cX) $ as claimed in the lemma.
\end{proof}
\begin{notation}[The orbifold Todd class of $ T\cX $]\label{notation_todd}
	Let $ \cX = [X/G] $ be a connected separated smooth quotient DM stack.
	We have
	\begin{equation*}
	\orbtd(T\cX) = \frac{\td(TI\cX)}{e^\rho(NI\cX)} = \left(\td(T\cX),\ \bigoplus_{0 \neq i \in I} \frac{\td(T\cX_i)}{e^{\rho_{g_i}}(N\cX_i)}\right)
	\end{equation*}
	in $ A(I\cX)_\CC^\times $, where
	\begin{equation*}
	\td(T\cX) = \frac{\td(TX)}{\td(\cO_X \otimes \fg)}
	\end{equation*}
	is the Todd class of the tangent complex of $ \cX $ in $ A(\cX)_\QQ^\times $. To simplify notations, we write
	\begin{equation*}
	\orbtd_\cX = \frac{\td_{I\cX}}{e^\rho_{I\cX}} = \left( \td_\cX, \bigoplus_{0 \neq i \in I} \frac{\td_{\cX_i}}{e^\rho_{\cX_i}} \right), \vspace{5pt}
	\end{equation*}
	where $ \td_{I\cX} = \td(TI\cX), \ e^\rho_{I\cX} = e^\rho(NI\cX), \ \td_{\cX_i} = \td(T\cX_i) $, and $ e^\rho_{\cX_i} = e^{\rho_{g_i}}(N\cX_i) $ for each $ 0 \neq i \in I. $
	If $ G $ is a finite group, then $ \fg $ vanishes and hence we have
	\begin{equation*}\label{eq_orbtd_finite_G}
	\orbtd_\cX = \left(\td_X, \bigoplus_{0 \neq i \in I} \frac{\td_{X_i}}{e^{\rho_{g_i}}(NX_i)}\right) = \left(\td_X, \bigoplus_{0 \neq i \in I} \frac{\td_{X_i}}{\prod_{\lambda, j} \left(1 - \lambda^{-1} e^{-h_{i,\lambda,j}}\right)}\right),
	\end{equation*}
	where the second identity is from the $ G_i $-equivariant decomposition
	\begin{equation*}
	NX_i = \bigoplus_{\lambda \in \eig(\tau_i(g_i))} N_{i,\lambda}
	\end{equation*}
	of the normal bundle $ NX_i $ for each $ 0 \neq i \in I $, where $ \tau_i $ is the $ G_i $-equivariant structure on $ NX_i $ induced from the $ G $-equivariant tangent bundle $ (TX, \tau) $ on $ X, $ $ \eig(\tau_i(g_i)) $ contains nontrivial $ \ord(g_i) $-th roots of unity, and $ h_{i,\lambda,1}, \dots, h_{i,\lambda,\rk(N_{i,\lambda})} $ are the Chern roots of each eigenbundle $ N_{i,\lambda} $ for each $ 0 \neq i \in I $.
\end{notation}
\vskip 2em

\section{Orbifold Mukai Pairing and Orbifold HRR Formula}\label{Orbifold_Mukai}
In this section we will define an orbifold Mukai pairing and derive an HRR formula for connected proper smooth quotient DM stacks. 
\vskip 4pt 
We begin with the definition of a proper algebraic stack.
\begin{defn}
	An algebraic stack $ \cX $ is proper if $ \cX $ is separated and of finite type, and if there exists a proper scheme $ Y $ with a surjective morphism $ Y \to \cX. $
\end{defn}
\begin{defn}\label{defn_Euler_char}
	Let $ \cX $ be a proper quotient DM stack. The \tb{orbifold Euler characteristic} $ \chi(\cX, \ {\cdot}\ ): K(\cX) \to \ZZ $ is defined by
	\begin{equation*}
	\chi(\cX, \cE) = \sum_i (-1)^i \dim H^i(\cX,\cE)
	\end{equation*}
	for a coherent sheaf $ \cE $ on $ \cX $ and extended linearly. 
\end{defn}
\begin{notation}
	Let $ \cX $ be a connected proper quotient DM stack of dimension $ n $ with its structure morphism $ f: \cX \to \pt $. We denote the pushforward 
	$$ f_*: A(\cX)_\CC \to A(\pt)_\CC = \CC $$ 
	by the integral symbol $\int_\cX $, i.e., for an element $ cv \in A(\cX)_\CC $ with $ c \in \CC $ and $ v \in A(\cX) $, we have
	\begin{equation*}
	\int_\cX cv = c\deg (v_n),
	\end{equation*}
	where $ v_n $ is the degree $ n $ component of $ v. $
	We also have a linear map
	\begin{equation*}
	\int_{I\cX} : A(I\cX)_\CC \to \CC	
	\end{equation*}
	defined by $ \int_{I\cX} = \sum_{i \in I} \int_{\cX_i} $. 
\end{notation}
Now we can state the orbifold Hirzebruch-Riemann-Roch (HRR) theorem.
\begin{thm}[Orbifold HRR]\label{thm_HRR}
	Let $ \cX $ be a connected proper smooth quotient DM stack. Then for all $ x \in K(\cX), $ we have
	\begin{equation}\label{eq_HRR}
	\chi(\cX, x) = \int_{I\cX} \orbch(x) \orbtd_\cX = \int_{I\cX} \orbch(x) \frac{\td_{I\cX}}{e^\rho_{I\cX}}.
	\end{equation}
\end{thm}
\begin{proof}
	Take an element $ x $ in $ K(\cX). $ By the natural isomorphism  $ K(\cX) \cong K^0(\cX) $, we can write $ x = [\cV] - [\cW] $ for two vector bundles $ \cV $ and $ \cW $ on $ \cX $. Formula (\ref{eq_HRR}) follows from Theorem 4.19 in \cite{edidin2013riemann} for vector bundles on $ \cX $ and the linearities of the orbifold Euler characteristic and the orbifold Chern character.
\end{proof}
\begin{rmk}[The explicit orbifold HRR theorem]
	Let $ I $ be the finite inertia index set of a connected proper smooth quotient DM stack $ \cX $ in Notations \ref{notation_quotient_stack} and \ref{notation_todd}. For a vector bundle $ \cV $ on $ \cX $, the orbifold HRR theorem becomes
	\begin{equation*}
	\chi(\cX, \cV) = \int_{\cX} \ch(\cV) \td_\cX + \sum_{0 \neq i \in I} \int_{\cX_i} \ch(\cV_i) \frac{\td_{\cX_i}}{e^{\rho}_{\cX_i}},
	\end{equation*}
	where $ \cV_i $ is the restriction of the vector bundle $ \cV $ on the substack $ \cX_i $ of $ \cX $ for each $ 0 \neq i \in I. $
\end{rmk}
\begin{defn}\label{defn_Euler_pairing}
	Let $ \cX = [X/G] $ be a proper quotient stack. The \tb{orbifold Euler pairing} 
	$$ \chi: K(\cX) \times K(\cX) \to \ZZ $$ 
	is defined by
	\begin{equation*}
	\chi(\cE,\cF) = \sum_i (-1)^i \dim \Ext^i(\cE,\cF)
	\end{equation*}
	for two coherent sheaves $ \cE $ and $ \cF $ on $ \cX $ and extended bilinearly.
\end{defn}
\begin{lem}\label{lem_Euler_pairing}
	Let $ \cX = [X/G] $ be a proper smooth quotient DM stack. Then for all $ x $ and $ y $ in $ K(\cX), $ we have
	\begin{equation*}
	\chi(x, y) = \chi(\cX, x^\vee y).
	\end{equation*}
	In particular, $ \chi(1, x) = \chi(\cX, x) $ for all $ x $ in $ K(\cX) $, where $ 1 = [\cO_\cX]. $
\end{lem}
\begin{proof}
	By the bilinearity of $ \chi, $ it suffices to consider $ x = [\cV] $ and $ y =[\cW] $ for vector bundles $ \cV = (V, \phi) $ and $ \cW = (W, \psi) $ on $ \cX. $ Then we have
	\begin{align*}
	\chi(x, y) & = \chi(\cV, \cW) = \sum_i (-1)^i \dim \Ext^i(V,W)^G \\
	& = \sum_i (-1)^i \dim H^i(X,V^\vee \otimes W)^G & \text{because $ V $ is locally free} \\
	& = \chi(\cX, \cV^\vee \otimes \cW) \\
	& = \chi(\cX, x^\vee y).
	\end{align*}
\end{proof}

Next we will define the orbifold Mukai vector and the orbifold Mukai pairing for quotient stacks.

Let $ \cX $ be a connected separated smooth quotient DM stack. We first define the square root of the Todd class of the inertia stack $ I\cX. $
\begin{defn}
	Let $ \cX $ be a connected separated smooth quotient DM stack with an inertia index set $ I $ in Notations \ref{notation_quotient_stack} and \ref{notation_todd}.
	Take an index $ i \in I. $ We have $ \td_{\cX_i} = 1+v_i \in A(\cX_i)_\CC^\times $ for some $ v_i \in A^{\geq 1}(\cX_i). $ If $ d_i = \dim \cX_i $, then the square root of $ \td_{\cX_i} $
	$$ \sqrt{\td_{\cX_i}} = \sum_{k=0}^{d_i} \binom{1/2}{k} v_i^k = 1 + \frac{v_i}{2} - \frac{v_i^2}{8} + \cdots + \binom{1/2}{d_i} v_i^{d_i} $$
	is well defined in $ A(\cX_i)_\CC^\times. $ We define the square root of $ \td_{I\cX} $ by
	\begin{equation*}
		\sqrt{\td_{I\cX}} = \bigoplus_{i \in I} \sqrt{\td_{\cX_i}}.
	\end{equation*}
\end{defn}

\begin{defn}\label{defn_orbv}
	Let $ \cX $ be a connected separated smooth quotient DM stack. The \tb{orbifold Mukai vector} map  $ \orbv: K(\cX) \to A(I\cX)_\CC $
	is defined by
	\begin{equation*}
	\orbv(x) = \orbch(x) \sqrt{\td_{I\cX}}
	\end{equation*}
	for all $ x $ in $ K(\cX). $
\end{defn}
\begin{rmk}
	For an element $ x $ in $ K(\cX), $ if $ \orbv(x) = (v,(v_i)) \in A(\cX)_\QQ \oplus A(I_t\cX)_\CC $, then $ v $ is called the \tb{Mukai vector} of $ x $ and each $ v_i $ for $ 0 \neq i \in I $ is called a \tb{twisted Mukai vector} of $ x $.
\end{rmk}
\begin{rmk}
	Our definition of the orbifold Mukai vector map is different from the one in \cite{popa2017derived}, where
	\begin{equation*}
	\orbv(x) = \orbch(x) \sqrt{\orbtd_{\cX}} = \orbch(x) \sqrt{\td_{I\cX}/e^\rho_{I\cX}}
	\end{equation*}
	in the cohomology $ H^*(I\cX)_\CC. $ We made a different (and arguably better) choice because taking the square root of $ e^\rho_{I\cX} $ is an unnecessary computation and it would also make the orbifold Mukai pairing more complicated. We also use the complex Chow ring $ A(I\cX)_\CC $ rather than $ H^*(I\cX)_\CC $ for the target of the orbifold Mukai vector map.
\end{rmk}
The next result follows immediately from Definition \ref{defn_orbv}.
\begin{prop}
	Let $ \cX $ be a connected separated smooth quotient DM stack. The orbifold Mukai vector map $ \orbv: K(\cX) \to A(I\cX)_\CC $ satisfies
	\begin{equation*}
	\orbv(x+y) = \orbv(x) + \orbv(y) \quad \text{and} \quad \orbv(xy) = \orbv(x) \orbch(y)
	\end{equation*}
	for all $ x $ and $ y $ in $ K(\cX). $ 
\end{prop}
There is an involution on the complex Chow ring $ A(I\cX)_\CC $.
\begin{defn}\label{defn_involution}
	Let $ \cX $ be a connected separated smooth quotient DM stack. We define an involution 
	$$ (\ \cdot\ )^\vee: A(I\cX)_\CC \to A(I\cX)_\CC $$ 
	by defining one on each component of $ A(I\cX)_\CC $. Take $ i \in I. $ An element $ u \in A(\cX_i)_\CC $ is a finite sum:
	\begin{equation*}
	u = \sum_{j,k} a_{jk} u_{jk},
	\end{equation*}
	where $ a_{jk} \in \CC $, and $ u_{jk} \in A^j(\cX_i) $ is the class of an irreducible closed substack of $ \cX_i $ of codimension $ j $. The involution on $ A(\cX_i)_\CC $ is defined by
	\begin{equation*}
	u^\vee = \sum_{j, k} (-1)^j \overline{a}_{jk} u_{jk},	
	\end{equation*}
	where $ \overline{a}_{jk} $ denotes the complex conjugate of $ a_{jk} $. For every element $ \orbv = \bigoplus_{i \in I} v_i $ in $ A(I\cX)_\CC, $ define
	\begin{align*}
	\orbv^\vee = \bigoplus_{i \in I} v_i^\vee.
	\end{align*}
\end{defn}
\begin{rmk}
	It is straightforward to check that $ (\ \cdot \ )^\vee: A(I\cX)_\CC \to A(I\cX)_\CC $ is a ring automorphism which also commutes with $ \sqrt{(\ \cdot \ )} $ when this is defined.
\end{rmk}
\begin{prop}\label{lem_orbch_orbv}
	Let $ \cX $ be a connected separated smooth quotient DM stack. Then for all $ x $ and $ y$ in $ K(\cX), $ we have
	\begin{equation*}
	\orbch(x^\vee) = \orbch(x)^\vee \quad \text{and} \quad \orbv(x^\vee) = \orbv(x)^\vee \sqrt{\td_{I\cX}/\td_{I\cX}^\vee}
	\end{equation*}
	in $ A(I\cX)_\CC $.
\end{prop}
\begin{proof}
	We first prove the first identity. By the linearity of $ \orbch $ and the involutions on $ K(\cX) $ and $ A(I\cX)_\CC, $ it suffices to check the case $ x = [\cV] $ is the class of a vector bundle $ \cV $ on $ \cX. $ We can write
	$$ \orbch\left(\cV^\vee\right) = \left( \ch(\cV^\vee),\ \bigoplus_{0 \neq i \in I} \ch^{\rho_{g_i}}\left(\cV_i^\vee\right) \right). $$ 
	We have $ \ch(\cV^\vee) = \ch(\cV)^\vee $ by the property of the ordinary Chern character. For each $ 0 \neq i \in I, $ we have
	\begin{align*}
	\ch^{\rho_{g_i}} \left(\cV_i^\vee\right) & = \ch \left(\sum_{\lambda \in \eig(\phi_i(g_i))} \lambda^{-1} \, [\cV_{i,\lambda}]^\vee\right) \\
	& = \sum_{\lambda \in \eig(\phi_i(g_i))} \sum_{j=0}^{\dim \cX_i} (-1)^j \lambda^{-1} \ch_j \left(\cV_{i,\lambda}\right) = \left(\ch^{\rho_{g_i}} (\cV_i)\right)^\vee.
	\end{align*}
	Now the second identity follows from two identities:
	\begin{equation*}
	\orbv (x^\vee) = \orbch(x^\vee) \sqrt{\td_{I\cX}} = \orbch(x)^\vee \sqrt{\td_{I\cX}} \quad \text{and}\quad \orbv (x)^\vee = \orbch(x)^\vee \sqrt{\td_{I\cX}}^\vee
	\end{equation*}
\end{proof}

Assume $ \cX $ is proper. Now we define an orbifold Mukai pairing on $ A(I\cX)_\CC $.
\begin{defn}\label{defn_orbv_pairing}
	Let $ \cX $ be a connected proper smooth quotient DM stack. The \tb{orbifold Mukai pairing} 
	$$ \inprod{\cdot \ {,} \ \cdot }_{I\cX}: A(I\cX)_\CC \times A(I\cX)_\CC \to \CC $$ 
	is defined by
	\begin{equation}\label{eq_orbv_pairing}
	\inprod{\orbv,\orbw}_{I\cX} = \int_{I\cX} \frac{\orbv^\vee \orbw}{e^\rho_{I\cX}} \sqrt{\frac{\td_{I\cX}}{\td_{I\cX}^\vee}}
	\end{equation}
	for all $ \orbv $ and $ \orbw $ in $ A(I\cX)_\CC. $
\end{defn}
\begin{rmk}
	For a pair of vectors $ \orbv = (v,(v_i)) $ and $ \orbw = (w, (w_i)) $ in $ A(I\cX)_\CC, $ the orbifold Mukai pairing can be written as
	\begin{equation*}
	\inprod{\orbv,\orbw}_{I\cX} = \inprod{v,w}_\cX + \sum_{0 \neq i \in I} \inprod{v_i,w_i}_{\cX_i},
	\end{equation*}
	where the \tb{Mukai pairing}
	\begin{equation*}
	\inprod{v,w}_\cX = \int_\cX v^\vee w \sqrt{\td_\cX/\td_\cX^\vee},
	\end{equation*}
	and the $ \boldsymbol{i} $\tb{-th twisted Mukai pairing}
	\begin{equation*}
	\inprod{v_i,w_i}_{\cX_i} = \int_{\cX_i} \frac{v_i^\vee w_i}{e^\rho_{\cX_i}} \sqrt{\td_{\cX_i}/\td_{\cX_i}^\vee}
	\end{equation*}
	for each $ 0 \neq i \in I. $
\end{rmk}
The following is immediate from the definition of the orbifold Mukai pairing.
\begin{prop}
	Let $ \cX $ be a connected proper smooth quotient DM stack. The orbifold Mukai pairing $ \inprod{\cdot \ {,} \ \cdot }_{I\cX}: A(I\cX)_\CC \times A(I\cX)_\CC \to \CC $ is sesquilinear, i.e.,
	\begin{equation*}
	\inprod{\orbv + \orbs, \orbw + \orbt}_{I\cX} = \inprod{\orbv, \orbw}_{I\cX} + \inprod{\orbv, \orbt}_{I\cX} + \inprod{\orbs, \orbw}_{I\cX} + \inprod{\orbs, \orbt}_{I\cX}
	\end{equation*}
	and
	\begin{equation*}
	\inprod{a\orbv, b\orbw}_{I\cX} = \cj{a}b \inprod{\orbv, \orbw}_{I\cX}
	\end{equation*}
	for all $ \orbv, \orbs, \orbw, \orbt $ in $ A(I\cX)_\CC $ and all $ a,b $ in $ \CC. $
\end{prop}
Now we can state the orbifold HRR formula in terms of the orbifold Euler pairing and the orbifold Mukai pairing.
\begin{thm}[Orbifold HRR formula]
	Let $ \cX $ be a connected proper smooth quotient DM stack. Then for all $ x $ and $ y $ in $ K(\cX), $ we have
	\begin{equation}\label{eq_HRR2}
	\chi(x, y) = \inprod{\orbv(x),\orbv(y)}_{I\cX}.
	\end{equation}
\end{thm}
\begin{proof}
	Take any $ x $ and $ y $ in $ K(\cX). $ Formula (\ref{eq_HRR2}) is a result of the orbifold HRR theorem and the definition of the orbifold Mukai pairing:
	\begin{align*}
	\chi(x,y) & = \chi(\cX, x^\vee y) & \text{by Lemma \ref{lem_Euler_pairing}} \\
	& = \int_{I\cX} \orbch(x^\vee y) \orbtd_\cX & \text{by the orbifold HRR theorem \ref{thm_HRR}} \\
	& = \int_{I\cX} \orbch(x^\vee) \orbch(y) \orbtd_\cX & \text{since $ \orbch $ is a ring map by Proposition \ref{prop_orbch_ring_map}} \\
	& = \int_{I\cX} \frac{\orbv(x^\vee) \orbv(y)}{e^\rho_{I\cX}} & \text{by Definition \ref{defn_orbv}} \\
	& = \int_{I\cX} \frac{\orbv(x)^\vee \orbv(y)}{e^\rho_{I\cX}} \sqrt{\td_{I\cX}/\td_{I\cX}^\vee} & \text{by Lemma \ref{lem_orbch_orbv}} \\
	& = \inprod{\orbv(x),\orbv(y)}_{I\cX} & \text{by Definition \ref{defn_orbv_pairing}.}
	\end{align*}
\end{proof}
The orbifold Euler pairing on $ K(\cX) $ can be extended to a sesquilinear form on $ K(\cX)_\CC $. Since $ \sqrt{\td_{I\cX}} $ is a unit in $ A(I\cX)_\CC, $ Proposition \ref{prop_orbch_isom} implies the following.
\begin{prop}
	Let $ \cX $ be a connected proper smooth quotient DM stack. The orbifold Mukai vector map $ \orbv: K(\cX)_\CC \to A(I\cX)_\CC $ is a linear isometry
	\begin{equation*}
	\orbv: \left(K(\cX)_\CC, \chi\right) \xrightarrow{\simeq} \left(A(I\cX)_\CC, \inprod{\cdot \ {,} \ \cdot }_{I\cX}\right).
	\end{equation*}
\end{prop}
Note that formula (\ref{eq_HRR2}) reduces to (\ref{eq_HRR}) when $ x = [\cO_\cX] = 1.$ Let's apply the orbifold HRR formula for $ \cX = BG $ for a finite group $ G. $
\begin{ex}
	Consider $ \cX = BG $ for a finite group $ G $ with conjugacy classes $ [g_0], [g_1], \dots, [g_{n-1}] $. Hence 
	$$ K(BG) = K^0(BG) = R(G). $$ 
	The orbifold Euler characteristic of a sheaf $ \cV = (V, \phi) $ on $ BG $ is
	\begin{equation*}
	\chi(BG, \phi) = \dim V^G,
	\end{equation*}
	and the orbifold Euler pairing for sheaves $ \cV = (V, \phi) $ and $ \cW = (W, \psi) $ on $ BG $ is
	\begin{equation*}
	\chi(\phi, \psi) = \dim \Hom(V,W)^G.
	\end{equation*}
	The inertia stack of $ BG $ is
	\begin{equation*}
	IBG \simeq \coprod_{i = 0}^{n-1} \{g_i\} \times Z_{g_i}, \quad \text{and} \quad A(IBG)_\CC \cong \CC^n.
	\end{equation*}
	The orbifold Mukai vector and orbifold Chern character coincide to be
	\begin{equation*}
	\orbv = \orbch: R(G) \to \CC^n, \quad [V, \phi] \mapsto \chi_\phi,
	\end{equation*}
	i.e., the character map for representations of $ G. $ The involution on $ R(G) $ is taking dual representations of $ G, $ and the involution on $ A(IBG)_\CC \cong \CC^n $ is complex conjugation. The orbifold Mukai pairing is an inner product on $ \CC^n $ weighted by the centralizers $ Z_{g_0}, \dots, Z_{g_{n-1}}, $ since formula (\ref{eq_orbv_pairing}) reduces to
	\begin{equation*}
	\inprod{a,b}_{IBG} = \int_{IBG} \overline{a} b = \sum_{i=0}^{n-1} \int_{BZ_{g_i}} \overline{a}_i b_i = \sum_{i=0}^{n-1} \frac{1}{|Z_{g_i}|} \overline{a}_i b_i
	\end{equation*}
	for all $ a = (a_0, \dots, a_{n-1}) $ and $ b = (b_0, \dots, b_{n-1}) $ in $ \CC^n. $
	For a coherent sheaf $ \cV = (V,\phi) $ on $ BG $, the orbifold HRR theorem reads
	\begin{align*}
	\chi(BG, \phi) = \dim V^G = \int_{IBG} \orbch(\cV) = \sum_{i=0}^{n-1} \int_{BZ_{g_i}} \chi_\phi(g_i) = \sum_{i=0}^{n-1} \frac{1}{|Z_{g_i}|} \chi_\phi(g_i) = \frac{1}{|G|} \sum_{g \in G} \chi_\phi(g),
	\end{align*}
	where the last equality holds because 
	$$ |[g]| \times |Z_{g}| = |G| $$ 
	for all $ g $ in $ G. $
	For a pair of coherent sheaves $ \cV = (V, \phi) $ and $ \cW = (W, \psi) $ on $ BG $, the orbifold HRR formula reads
	\begin{align*}
	\chi(\phi, \psi) = & \dim \Hom(V,W)^G = \inprod{\orbv(\cV),\orbv(\cW)}_{IBG} = \inprod{\chi_\phi,\chi_\psi}_W \\
	& = \sum_{i=0}^{n-1} \frac{1}{|Z_{g_i}|} \overline{\chi_\phi(g_i)}\chi_\psi(g_i) = \frac{1}{|G|} \sum_{g \in G} \chi_\phi(g^{-1}) \chi_\psi(g),
	\end{align*}
	which is a standard result in representation theory. When $ \phi $ is an irreducible representation of $ G, $ then $ \chi(\phi, \psi) $ counts the copies of $ \phi $ in $ \psi. $
\end{ex}
\begin{ex}[The orbifold HRR formula for $ B\mu_n $ is Parseval's theorem.]\label{ex_descrete_Parseval}
	Consider $ \cX = B\mu_n. $ The orbifold Euler pairing on $ K(B\mu_n) \cong \ZZ[x]/(x^n-1) $ is given by
	\begin{equation*}
	\inprod{x^i, x^j} = \delta_{ij}
	\end{equation*}
	for $ i,j = 0, \dots, n-1 $ and extended linearly. The orbifold Mukai pairing on
	$$ A(IB\mu_n) \cong \CC^n $$ 
	is the weighted inner product
	\begin{equation*}
	\inprod{a,b}_W = \frac{1}{n} \sum_{i=0}^{n-1} \cj{a}_i b_i
	\end{equation*}
	in $ \CC $ for all $ a = (a_0, \dots, a_{n-1}) $ and $ b = (b_0, \dots, b_{n-1}) $ in $ \CC^n. $
	Recall that the orbifold Mukai vector is the inverse discrete Fourier transform
	\begin{equation*}
	\orbv: \CC[x]/(x^n-1) \cong \CC^n \to \CC^n, \quad f \mapsto \check{f},
	\end{equation*}
	and its inverse is the discrete Fourier transform $ f \mapsto \hat{f}. $
	The orbifold HRR formula (\ref{eq_HRR2}) now reads
	\begin{equation*}
	\inprod{f, g} = \inprod{\check{f},\check{g}}_W
	\end{equation*}
	for all $ f, g $ in $ \ZZ[x]/(x^n-1). $ Extending the orbifold Euler pairing to $ K(B\mu_n)_\CC = \CC[x]/(x^n-1) $ sesquilinearly, we have
	\begin{equation*}
	\inprod{a, b}_W = \inprod{\hat{a},\hat{b}}
	\end{equation*}
	in $ \CC $ for all $ a, b $ in $ \CC^n, $
	which is Parseval's theorem for the discrete Fourier transform.
\end{ex}

%\bibliographystyle{alpha}
%\bibliography{ref}

\begin{thebibliography}{Kaw79}
	
	\bibitem[Edi13]{edidin2013riemann}
	Dan Edidin.
	\newblock Riemann-{R}och for {D}eligne-{M}umford stacks.
	\newblock In {\em A celebration of algebraic geometry}, volume~18 of {\em Clay
		Math. Proc.}, pages 241--266. Amer. Math. Soc., Providence, RI, 2013.
	
	\bibitem[EG98]{edidin1998equivariant}
	Dan Edidin and William Graham.
	\newblock Equivariant intersection theory (with an {A}ppendix by {A}ngelo
	{V}istoli: the {C}how ring of {M}2).
	\newblock {\em Invent. Math.}, 131(3):595--634, 1998.
	
	\bibitem[EG00]{edidin2000riemann}
	Dan Edidin and William Graham.
	\newblock Riemann-{R}och for equivariant {C}how groups.
	\newblock {\em Duke Math. J.}, 102(3):567--594, 2000.
	
	\bibitem[EG05]{edidin2005nonabelian}
	Dan Edidin and William Graham.
	\newblock Nonabelian localization in equivariant {$K$}-theory and
	{R}iemann-{R}och for quotients.
	\newblock {\em Adv. Math.}, 198(2):547--582, 2005.
	
	\bibitem[FL85]{fulton1985riemann}
	William Fulton and Serge Lang.
	\newblock {\em Riemann-{R}och {A}lgebra}, volume 277 of {\em Grundlehren der
		mathematischen Wissenschaften}.
	\newblock Springer New York, NY, 1985.
	
	\bibitem[Ful98]{fulton1998intersection}
	William Fulton.
	\newblock {\em Intersection {T}heory}, volume~2 of {\em Ergebnisse der
		Mathematik und ihrer Grenzgebiete. 3. Folge. A Series of Modern Surveys in
		Mathematics}.
	\newblock Springer New York, NY, second edition, 1998.
	
	\bibitem[Iri09]{iritani2009integral}
	Hiroshi Iritani.
	\newblock An integral structure in quantum cohomology and mirror symmetry for
	toric orbifolds.
	\newblock {\em Adv. Math.}, 222(3):1016--1079, 2009.
	
	\bibitem[Kaw79]{kawasaki1979riemann}
	Tetsuro Kawasaki.
	\newblock The {R}iemann-{R}och theorem for complex {$V$}-manifolds.
	\newblock {\em Osaka Math. J.}, 16(1):151--159, 1979.
	
	\bibitem[Pop17]{popa2017derived}
	Mihnea Popa.
	\newblock Derived equivalences of smooth stacks and orbifold {H}odge numbers.
	\newblock In {\em Higher {D}imensional {A}lgebraic {G}eometry: in honour of
		{P}rofessor {Y}ujiro {K}awamata's sixtieth birthday}, volume~74 of {\em Adv.
		Stud. Pure Math.}, pages 357--380. Math. Soc. Japan, Tokyo, 2017.
	
	\bibitem[To{\"e}99]{toen1999theoremes}
	Bertrand To{\"e}n.
	\newblock Th\'{e}or\`emes de {R}iemann-{R}och pour les champs de
	{D}eligne-{M}umford.
	\newblock {\em $K$-Theory}, 18(1):33--76, 1999.
	
	\bibitem[To{\"e}00]{toen2000motives}
	Bertrand To{\"e}n.
	\newblock On motives for {D}eligne-{M}umford stacks.
	\newblock {\em Internat. Math. Res. Notices}, 2000(17):909--928, 2000.
	
	\bibitem[Tot04]{totaro2004resolution}
	Burt Totaro.
	\newblock The resolution property for schemes and stacks.
	\newblock {\em J. Reine Angew. Math.}, 577:1--22, 2004.
	
	\bibitem[Tse10]{tseng2010orbifold}
	Hsian-Hua Tseng.
	\newblock Orbifold quantum {R}iemann-{R}och, {L}efschetz and {S}erre.
	\newblock {\em Geom. Topol.}, 14(1):1--81, 2010.
	
	\bibitem[Yau10]{yau2010lambda}
	Donald Yau.
	\newblock {\em Lambda-rings}.
	\newblock World Scientific Publishing Co. Pte. Ltd., Hackensack, NJ, 2010.
	
\end{thebibliography}

\end{document}